\documentclass[a4paper,10pt]{amsart}
\usepackage[german,english]{babel}
\usepackage[utf8x]{inputenc}
\usepackage{ucs}
\usepackage{amsmath}
\usepackage{amsfonts}
\usepackage{amssymb}
\usepackage{amsthm}
\usepackage{fontenc}
\usepackage{graphicx}
\usepackage{url}
\usepackage{color}
\usepackage[T1]{fontenc}
\usepackage{a4wide}

\newtheorem{theorem}{Theorem}

\newtheorem{lemma}{Lemma}

\newtheorem{definition}{Definition}

\newtheorem{remark}{Remark}

\title[A Posteriori Estimates For Time-Periodic Eddy Current Problems]{
A posteriori error estimation for the optimal control of time-periodic eddy current problems}

\author{M. Wolfmayr}

\address[M. Wolfmayr]{Jamk University of Applied Sciences, Finland, and University of Jyv\"{a}skyl\"{a}, Finland}

\email{monika.wolfmayr@jamk.fi}

\begin{document}

\begin{abstract}
This work presents the multiharmonic analysis and derivation of functional type a posteriori estimates of a distributed eddy current optimal control problem and its state equation in a time-periodic setting. The existence and uniqueness of the solution of a weak space-time variational formulation for the optimality system and the forward problem are proved by deriving inf-sup and sup-sup conditions. Using the inf-sup and sup-sup conditions, we derive guaranteed, sharp, and fully computable bounds of the approximation error for the optimal control problem and the forward problem in the functional type a posteriori estimation framework. We present here the first computational results on the derived estimates.
\end{abstract}
\maketitle

\section{Introduction}
\label{Sec1:Introduction}

We discuss time-periodic eddy current optimal control and corresponding boundary value problems. We derive a posteriori error estimates for the time-periodic boundary value problem and the optimality system of the optimal control problem.
Similar estimates have been derived previously for time-periodic parabolic problems in \cite{LangerRepinWolfmayr:2015, LangerRepinWolfmayr:2016}. However, functional type estimates for the time-periodic eddy-current problems are new and discussed here.
The problems occur in the application of Maxwell equations in electromagnetism. Eddy-current models for Maxwell equations haven been discussed in \cite{BuffaAmmariNedelec:2000}.

The problems are formulated in terms of their Fourier series expansions in time which is a natural framework due to the time-periodicity.
The Fourier coefficients can be then discretized by for instance the finite element method. This method has been applied previously to nonlinear eddy-current problems in \cite{YamadaBessho:1988} and later analyzed and efficient solvers presented in \cite{ BachingerKaltenbacherReitzinger:2002, BachingerLangerSchoeberl:2005, BachingerLangerSchoeberl:2006}. Domain decomposition methods have been applied in \cite{CopelandKolmbauerLanger:2010, CopelandLanger:2010}.

Optimal control problems are subject matter of many publications mentioning the books \cite{HinzePinnauUlbrichUlbrich:2009, Troeltzsch:2010, BorziSchulz:2012}.
The multiharmonic finite element method is a type of space-time method. Space-time methods for parabolic optimal control problems have been recently presented in \cite{LangerSteinbachTroeltzschYang:2021, LangerSteinbachTroeltzschYang:2021b, LangerSchafelner:2022}.

The multiharmonic analysis with inf-sup and sup-sup estimates regarding time-periodic parabolic optimal control problems have been discussed in \cite{LangerWolfmayr:2013, Wolfmayr:2014}. 
The conditions yield existence and uniqueness of a solution by applying Babu\v{s}ka-Aziz' theorem, see \cite{Babuska:1971, BabuskaAziz:1972}.
Robust preconditioning for the MINRES (minimal residual) method, which was introduced in \cite{PaigeSaunders:1975}, was presented in \cite{KollmannKolmbauer:2013, KollmannKolmbauerLangerWolfmayrZulehner:2013} for time-periodic and multiharmonic parabolic optimal control problems.

Functional type a posteriori estimates have been introduced in e.g. \cite{Repin:1999} and discussed for time-dependent problems in \cite{Repin:2002}.
These type of reliable methods have been discussed in various papers and books. We refer to the books \cite{NeittaanmaekiRepin:2004, Repin:2008, MaliNeittaanmaekiRepin:2013} and recently \cite{ PraetoriusRepinSauter:2021} discussing also eddy-current problems.
Functional a posteriori estimates have been applied to Maxwell type problems in \cite{Hannukainen:2008, AnjamMaliMuzalevskyNeittaanmaekiRepin:2009, NeittaanmakiRepin:2010, PaulyRepinRossi:2011, Anjam:2014}.
For optimal control problems, functional type a posteriori estimates have been discussed in \cite{GaevskayaHoppeRepin:2006, GaevskayaHoppeRepin:2007} together with the minimization of quadratic functionals with respect to parameters introduced from Young's inequality.

The focus of this work lies on eddy-current problems due to the range of applications regarding simulations of electromagnetic devices. In order to compute candidates for the approximations of the exact solution, we have discretized the problem by the multiharmonic finite element method. The multiharmonic method has been applied for time-discretization and the space-time dependent systems of linear equations corresponding to the Fourier coefficients have been discretized by the finite element method. 
In \cite{Yousept:2012}, a finite element analysis for the coefficients in a time-harmonic setting has been discussed.
We have used N\'{e}d\'{e}lec (edge) basis functions of lowest order for approximating the space curl-space, see \cite{Nedelec:1980, Nedelec:1986}.
For the numerical tests, we have used the \textit{Fast FEM assembly: edge elements toolbox} for computing the mass, stiffness matrices and the load vector, see \cite{RahmanValdman:2013, AnjamValdman:2015, Valdman:2023}.
Efficient solvers and preconditioners for time-periodic eddy-current optimal control problems have already been discussed in \cite{KolmbauerLanger:2012, Kolmbauer:2012:thesis, KolmbauerLanger:2013} or later in  \cite{AxelssonLiang:2019, AxelssonLukas:2019}.

Adaptive methods for eddy current problems have been discussed in e.g. \cite{ChenChenCuiZhang:2010}.
In \cite{XuYouseptZou:2020}, an adaptive edge element method for a quasilinear curl curl problem has been presented.
Hierarchical error estimators for eddy-current problems are subject matter of \cite{BeckHiptmairWohlmuth:1999} and
residual based error estimators of e.g. \cite{BeckHiptmairHoppeWohlmuth:2000, Schoeberl:2008} and recently in
\cite{BoffiGastaldiRodriguezSebestov:2017}.
Recent works on error estimators for harmonic eddy-current problems include \cite{CreuseLeMenachNicaise:2019}.
In \cite{ChaumontFreletErnVohralik:2022}, estimates by broken patchwise equilibration are presented.

The paper is organized as follows: 
section \ref{Sec2:OptimalControlProblem} presents the model problem, which is a minimization problem with respect to state and control and a time-periodic state equation, the forward problem. 
In section \ref{Sec3:Setting}, we present the Fourier space framework for the time-periodic problem, and derive the weak space-time variational formulations which form the basis for the inf-sup and sup-sup conditions as well as estimates in section \ref{Sec4:WSTVF}. 
The multiharmonic finite element discretization is presented in section \ref{Sec5:MFEM}. 
Finally, the a posteriori error estimates for the forward and optimal control problem are derived in sections \ref{Sec6:FunctionalAPosterioriEstimates:Forwardpde} and \ref{Sec7:FunctionalAPosterioriEstimates:OptiSys}, respectively.
The numerical results are discussed in section \ref{Sec8:NumericalResults} and conclusions and future outlook in section \ref{Sec9:ConclusionsOutlook}.

\section{Model problem}
\label{Sec2:OptimalControlProblem}

Let $\Omega \subset \mathbb{R}^3$ be a bounded Lipschitz domain, where
$\Gamma := \partial \Omega$ denotes the boundary.
We consider a time-periodic setting with a given time interval $(0,T)$
which equals also the time period.
We denote by $Q := \Omega \times (0,T)$ and $\Sigma := \Gamma \times (0,T)$
the full space-time domain, also called space-time cylinder,
and its outer surface area, respectively.
Let $\boldsymbol{y_d}$ be the given desired state.
The state and control functions $\boldsymbol{y}$ and 
$\boldsymbol{u}$ are subject of the minimization problem
\begin{align}
\label{equation:minfunc}
 \min_{\boldsymbol{y},\boldsymbol{u}} \mathcal{J}(\boldsymbol{y},\boldsymbol{u})
 =  \min_{\boldsymbol{y},\boldsymbol{u}}
 \frac{1}{2} \| \boldsymbol{y} - \boldsymbol{y_d} \|_{L^2(Q)}^2
 + \frac{\alpha}{2} \| \boldsymbol{u}\|_{L^2(Q)}^2 
\end{align}
and time-periodic state equation
\begin{align}
\label{equation:forwardpde:pde}
 \sigma \frac{\partial \boldsymbol{y}}{\partial t} + \textbf{curl}( \nu \, \textbf{curl} \, \boldsymbol{y})
	&= \boldsymbol{u} \hspace{1cm} &&\text{in } Q, \\
\label{equation:forwardpde:boundary}
 \boldsymbol{y} \times \boldsymbol{n} &= 0 \hspace{1cm} &&\text{on } \Sigma, \\
\label{equation:forwardpde:periodic}
 \boldsymbol{y}(0) &= \boldsymbol{y}(T) \hspace{1cm} &&\text{in } \overline{\Omega}, \\
 \label{equation:forwardpde:div}
 \text{div}(\sigma \boldsymbol{y}) &= 0 \hspace{1cm} &&\text{in } Q.
\end{align}
We incorporate the Coulomb gauging condition
\eqref{equation:forwardpde:div}
implicitly: the given desired state $\boldsymbol{y_d}$ is assumed to be weakly divergence-free, i.e.,
\begin{align}
\label{equation:divergencefree:yd}
 \int_\Omega \boldsymbol{y}_d \cdot \nabla v \, d\boldsymbol{x} = 0 \qquad \forall \, v \in H^1_0(\Omega).
\end{align}
The divergence-free property is passed over to the
unknown functions
$\boldsymbol{y}$ and $\boldsymbol{u}$ 
for eddy current optimal control problems
in a time-periodic setting.
The coefficients $\sigma = \sigma(\boldsymbol{x})$ and $\nu = \nu(\boldsymbol{x})$ denote the conductivity and reluctivity, respectively.
Altogether the problem forms a time-periodic distributed eddy current optimal control problem.
In this work,
we consider the linear version of the eddy current problem \eqref{equation:forwardpde:pde}-\eqref{equation:forwardpde:periodic}, where $\Omega$ is a conducting domain:
the reluctivity $\nu$ is independent of $|\textbf{curl} \, \boldsymbol{y}|$,
and $\sigma$ and $\nu$ are both strictly positive and uniformly bounded, i.e.,
\begin{align*}
 0 < \underline{\sigma} \leq \sigma(\boldsymbol{x}) \leq \overline{\sigma} \quad \mbox{ and } \quad
 0 < \underline{\nu} \leq \nu(\boldsymbol{x}) \leq \overline{\nu}, \quad \boldsymbol{x} \in \Omega.
\end{align*}
However, we could extend our analysis also to the case of bounded domains which consist of conducting and non-conducting domains, see \cite{KolmbauerLanger:2012}.
The desired state $\boldsymbol{y_d}$ is the given target that we try to reach via a suitable control $\boldsymbol{u}$.
The positive regularization parameter $\alpha$ provides a weighting of the cost of the control in the cost functional $\mathcal{J}(\cdot,\cdot)$ in \eqref{equation:minfunc}.

The solution of the optimal control problem \eqref{equation:minfunc}-\eqref{equation:forwardpde:periodic} is equivalent to the solution of its optimality system. The Lagrange functional of the minimization problem is given as
\begin{align}
\label{def:LagrangeFunc}
 \mathcal{L}(\boldsymbol{y},\boldsymbol{u},\boldsymbol{p}) := 
 \mathcal{J}(\boldsymbol{y},\boldsymbol{u}) 
 - \int_0^T \int_{\Omega} \big( \sigma \frac{\partial \boldsymbol{y}}{\partial t} + \textbf{curl}( \nu \, \textbf{curl} \, \boldsymbol{y})
	- \boldsymbol{u}\big) \boldsymbol{p} \, d\boldsymbol{x} dt,
\end{align}
where $\boldsymbol{p}$ denotes the so-called Lagrange multiplier or adjoint state.
The necessary optimality conditions
\begin{align}
 \label{equation:NecOptConds}
\begin{aligned}
 \nabla_{\boldsymbol{y}} \mathcal{L}(\boldsymbol{y},\boldsymbol{u},\boldsymbol{p}) = 0, \qquad
 \nabla_{\boldsymbol{u}} \mathcal{L}(\boldsymbol{y},\boldsymbol{u},\boldsymbol{p}) = 0,  \qquad
 \nabla_{\boldsymbol{p}} \mathcal{L}(\boldsymbol{y},\boldsymbol{u},\boldsymbol{p}) = 0,
\end{aligned}
\end{align}
characterize a stationary point of the Lagrange functional (\ref{def:LagrangeFunc}).
We do not impose any inequality constraints on the control.
Hence, we can use the second condition of (\ref{equation:NecOptConds}), i.e.,
$\boldsymbol{u}= - \alpha^{-1} \boldsymbol{p}$ in $Q$, 
to eliminate the control $\boldsymbol{u}$ from the optimality system.
This leads to the derivation of a reduced optimality system, which can be written in its classical formulation as
\begin{align}
\label{equation:KKTSysClassical:1}
 \sigma \frac{\partial \boldsymbol{y}}{\partial t} + \textbf{curl}( \nu \, \textbf{curl} \, \boldsymbol{y})
	&= - \alpha^{-1} \boldsymbol{p
	} \hspace{1cm} &&\text{in } Q, \\
\label{equation:KKTSysClassical:2}
 \boldsymbol{y} \times \boldsymbol{n} &= 0 \hspace{1cm} &&\text{on } \Sigma, \\
 \label{equation:KKTSysClassical:3}
 \boldsymbol{y}(0) &= \boldsymbol{y}(T) \hspace{1cm} &&\text{in } \overline{\Omega}, \\
 \label{equation:KKTSysClassical:4}
-  \sigma \frac{\partial \boldsymbol{p}}{\partial t} + \textbf{curl}( \nu \, \textbf{curl} \, \boldsymbol{p})
	&= \boldsymbol{y} - \boldsymbol{y_d} \hspace{1cm} &&\text{in } Q, \\
\label{equation:KKTSysClassical:5}
 \boldsymbol{p} \times \boldsymbol{n} &= 0 \hspace{1cm} &&\text{on } \Sigma, \\
 \label{equation:KKTSysClassical:6}
 \boldsymbol{p}(0) &= \boldsymbol{p}(T) \hspace{1cm} &&\text{in } \overline{\Omega}.
\end{align}
The weakly divergence-free property of $\boldsymbol{y}_d$ is passed on to the state $\boldsymbol{y}$ and the control $\boldsymbol{u}$, hence also to the adjoint state $\boldsymbol{p}$.

In this work, we derive a weak space-time variational formulation of the state equation \eqref{equation:forwardpde:pde}-\eqref{equation:forwardpde:div} and of the optimality system \eqref{equation:KKTSysClassical:1}-\eqref{equation:KKTSysClassical:6}. We then deduce a posteriori error estimates for the optimal control problem.
In the following section, we introduce Sobolev spaces of functions in the space-time domain $Q$. 
The notation is close to the one used by Ladyzhenskaya \cite{Ladyzhenskaya:1973} and Ladyzhenskaya et al. \cite{LadyzhenskayaSolonnikovUralceva:1968}.

\section{Fourier space setting}
\label{Sec3:Setting}

Let $\boldsymbol{L^2}(\Omega) := [L^2(\Omega)]^3$ and $\boldsymbol{L^2}(Q) := [L^2(Q)]^3$.
We define the spaces $\boldsymbol{H}^{\textbf{curl}}(\Omega) := \{ \boldsymbol{v} \in \boldsymbol{L^2}(\Omega) :
\textbf{curl} \, \boldsymbol{v} \in \boldsymbol{L^2}(\Omega) \}$ and
$\boldsymbol{H}^{\textbf{curl}}_0(\Omega) := \{ \boldsymbol{v} \in \boldsymbol{H}^{\textbf{curl}}(\Omega) :
\boldsymbol{v} \times \boldsymbol{n} = 0 \; \mbox{on} \;  \Gamma \}$.
Also, $H^1(\Omega) := \{ v \in L^2(\Omega) :  \nabla v \in \boldsymbol{L^2}(\Omega) \}$
and
$H^1_0(\Omega) := \{ v \in H^1(\Omega) :  v = 0 \; \mbox{on} \;  \Gamma \}$.

For the space-time domain $Q$, we define the function spaces 
$\boldsymbol{H}^{\textbf{curl},0}(Q) = \{ \boldsymbol{v} \in \boldsymbol{L^2}(Q) :
\textbf{curl}_{\boldsymbol{x}} \, \boldsymbol{v} \in \boldsymbol{L^2}(Q) \}	$
and
$\boldsymbol{H}^{\textbf{curl},1}(Q) = \{ \boldsymbol{v} \in \boldsymbol{L^2}(Q) :
\textbf{curl}_{\boldsymbol{x}} \, \boldsymbol{v} \in \boldsymbol{L^2}(Q),
\partial_t \boldsymbol{v} \in \boldsymbol{L^2}(Q) \}$.
Next, we define
$H^{1,0}(Q) = \{ v \in L^2(Q) :  \nabla_{\boldsymbol{x}} v \in \boldsymbol{L^2}(Q) \}$
and $H^{1,0}_0(Q)$, where the latter includes homogeneous Dirichlet boundary conditions.
The notations $\textbf{curl}_{\boldsymbol{x}}$, $\partial_t$ and $\nabla_{\boldsymbol{x}}$ are used for the weak spatial curl, the weak time derivative and the weak spatial gradient, respectively.

For ease of notation, we introduce the following simplifications: the index  $\boldsymbol{x}$ in $\textbf{curl}_{\boldsymbol{x}}$ and $\nabla_{\boldsymbol{x}}$ will be omitted, and all inner products and norms in $\boldsymbol{L^2}$ related to the whole space-time domain $Q$ will be denoted by $(\cdot,\cdot)$ and $\|\cdot\|$, respectively.
Inner product $(\cdot,\cdot)_{\Omega}$ and norm $\|\cdot\|_{\Omega}$ are related to $\boldsymbol{L^2}(\Omega)$,
inner product and norm $(\cdot,\cdot)_{1,\Omega}$ and $\|\cdot\|_{1,\Omega}$ to $\boldsymbol{H}^1(\Omega)$,
and $(\cdot,\cdot)_{\textbf{curl},\Omega}$ and $\|\cdot\|_{\textbf{curl},\Omega}$
to $\boldsymbol{H}^{\textbf{curl}}(\Omega)$.
All the inner product and norm symbols are used for both the scalar and the vector-valued case, since their use is obvious from the context and simplifies notation.
 
The space-time Sobolev spaces $\boldsymbol{H}^{\textbf{curl},0}(Q)$ and $\boldsymbol{H}^{\textbf{curl},1}(Q)$ are equipped with the norms
 \begin{align*}
  \|\boldsymbol{v}\|_{\boldsymbol{H}^{\textbf{curl},0}}^2 :=
  \|\boldsymbol{v} \|^2 +  \|\textbf{curl} \, \boldsymbol{v} \|^2
 \qquad
\text{and}
\qquad
  \|\boldsymbol{v}\|_{\boldsymbol{H}^{\textbf{curl},1}}^2 :=
  \|\boldsymbol{v} \|^2 +  \|\textbf{curl} \, \boldsymbol{v} \|^2
  +  \|\partial_t \boldsymbol{v} \|^2,
 \end{align*}
 respectively.
 Boundary and time-periodicity conditions are included by defining the following function spaces:
\begin{align*}
\boldsymbol{H}^{\textbf{curl},0}_0(Q) &:= \{ \boldsymbol{v} \in \boldsymbol{H}^{\textbf{curl},0}(Q): \boldsymbol{v} \times \boldsymbol{n} = 0 \; \mbox{on} \;  \Sigma \}, \\
\boldsymbol{H}^{\textbf{curl},1}_0(Q) &:= \{ \boldsymbol{v} \in \boldsymbol{H}^{\textbf{curl},1}(Q): \boldsymbol{v} \times \boldsymbol{n} = 0 \; \mbox{on} \;  \Sigma \}, \\
\boldsymbol{H}^{\textbf{curl},1}_{0,per}(Q) &:= \{ \boldsymbol{v} \in \boldsymbol{H}^{\textbf{curl},1}_0(Q): \boldsymbol{v}(\boldsymbol{x},0) = \boldsymbol{v}(\boldsymbol{x},T) \ \mbox{for almost all } \boldsymbol{x} \in \Omega \}.
\end{align*}
We introduce additionally the spaces
\begin{align*}
\boldsymbol{H}^{0,1}(Q) &:= \{ \boldsymbol{v} \in \boldsymbol{L^2}(Q): \partial_t \boldsymbol{v} \in \boldsymbol{L^2}(Q) \}, \\
\boldsymbol{H}^{0,1}_{per}(Q) &:= \{ \boldsymbol{v} \in \boldsymbol{H}^{0,1}(Q): \boldsymbol{v}(\boldsymbol{x},0) = \boldsymbol{v}(\boldsymbol{x},T) \ \mbox{for almost all } \boldsymbol{x} \in \Omega \}.
\end{align*}
All $L^2(Q)$ functions provide the feasible representation as Fourier series expansion.
The real Fourier series expansion in time for $\boldsymbol{v} \in \boldsymbol{L^2}(Q)$ is given by
\begin{align*}
\boldsymbol{v}(\boldsymbol{x},t) = \boldsymbol{v}_0^c(\boldsymbol{x}) + \sum_{k=1}^{\infty} \left(\boldsymbol{v}_k^c(\boldsymbol{x}) \cos(k \omega t) + \boldsymbol{v}_k^s(\boldsymbol{x}) \sin(k \omega t)\right)
\end{align*}
with the cosine and sine Fourier coefficients 
\begin{align*}
 \boldsymbol{v}_k^c(\boldsymbol{x}) = \frac{2}{T} \int_0^T \boldsymbol{v}(\boldsymbol{x},t) \cos(k \omega t)\,dt, \quad
 \boldsymbol{v}_k^s(\boldsymbol{x}) = \frac{2}{T} \int_0^T \boldsymbol{v}(\boldsymbol{x},t) \sin(k \omega t)\,dt
\end{align*}
and
\begin{align*}
 \boldsymbol{v}_k^c(\boldsymbol{x}) = \frac{1}{T} \int_0^T \boldsymbol{v}(\boldsymbol{x},t)\,dt.
\end{align*}
The periodicity is $T$ and the corresponding frequency $\omega = 2 \pi /T$.
We use the following notation:
\begin{align*}
\boldsymbol{v}_k = (\boldsymbol{v}_k^c,\boldsymbol{v}_k^s)^T, \quad
\boldsymbol{v}_k^{\perp} = (-\boldsymbol{v}_k^s,\boldsymbol{v}_k^c)^T, \quad
\textbf{curl} \, \boldsymbol{v}_k = (\textbf{curl} \, \boldsymbol{v}_k^c, \textbf{curl} \, \boldsymbol{v}_k^s)^T.
\end{align*}
The relation $\|\boldsymbol{v}_k^\perp\|^2_{\Omega} = \|\boldsymbol{v}_k\|^2_{\Omega}$ holds.
We define a perpendicular Fourier series as introduced in \cite{Wolfmayr:2014}:
\begin{align*}
\begin{aligned}
\boldsymbol{v}^{\perp}(\boldsymbol{x},t) = \sum_{k=1}^{\infty} \left(- \boldsymbol{v}_k^c(\boldsymbol{x}) \sin(k \omega t) + \boldsymbol{v}_k^s(\boldsymbol{x}) \cos(k \omega t)\right).
\end{aligned}
\end{align*}
We define the following $\sigma$-weighted inner products:
\begin{align}
\label{def:weighted-time-product}
 \big(\sigma 
 \partial_t \boldsymbol{y}, \boldsymbol{v}\big) :=
 \int_{Q} \sigma \partial_t \boldsymbol{y} \cdot
 \boldsymbol{v} \, d\boldsymbol{x}\,dt
= \frac{T}{2} \sum_{k=1}^{\infty} k \omega (\sigma 
 \boldsymbol{y}_k,\boldsymbol{v}_k)_{\Omega}
\end{align}
and
\begin{align}
\label{def:weighted-time-product:perp}
 \big(\sigma 
 \partial_t \boldsymbol{y}, \boldsymbol{v}^{\perp}\big) :=
 \int_{Q} \sigma \partial_t \boldsymbol{y} \cdot
\boldsymbol{v}^{\perp} \, d\boldsymbol{x}\,dt
= \frac{T}{2} \sum_{k=1}^{\infty} k \omega (\sigma 
 \boldsymbol{y}_k,\boldsymbol{v}_k^{\perp})_{\Omega}.
\end{align}
\begin{definition}
\label{definition:HalfSpaces}
By introducing the norm
$\big\| \partial^{1/2}_t \boldsymbol{v} \big\|^2 := \frac{T}{2} \sum_{k=1}^{\infty} k \omega \|\boldsymbol{v}_k\|_{\Omega}^2$
in Fourier space, we define the spaces
$\boldsymbol{H}^{0,\frac{1}{2}}_{per}(Q) :=  \{ \boldsymbol{v} \in \boldsymbol{L^2}(Q): \big\| \partial^{1/2}_t \boldsymbol{v} \big\| < \infty \}$,
$\boldsymbol{H}^{\emph{\textbf{curl}},\frac{1}{2}}_{per}(Q) := \{ \boldsymbol{v} \in \boldsymbol{H}^{\emph{\textbf{curl}},0}(Q): \big\| \partial^{1/2}_t \boldsymbol{v} \big\| < \infty \}$ and 
$\boldsymbol{H}^{\emph{\textbf{curl}},\frac{1}{2}}_{0,per}(Q) := \{ \boldsymbol{v} \in \boldsymbol{H}^{\emph{\textbf{curl}},\frac{1}{2}}_{per}(Q): \boldsymbol{v} \times \boldsymbol{n} = 0 \mbox{ on } \Sigma \}$ with
$|\boldsymbol{v}|_{\boldsymbol{H}^{0,\frac{1}{2}}_{per}} = \big\| \partial^{1/2}_t \boldsymbol{v} \big\|$.
We define the $\sigma$-weighted inner products as follows
\begin{align}
\label{def:weighted-half-time-product}
 \big(\sigma 
 \partial^{1/2}_t \boldsymbol{y}, \partial^{1/2}_t \boldsymbol{v}\big) :=
 \int_{Q} \sigma \partial_t^{1/2} \boldsymbol{y} \cdot
\partial_t^{1/2} \boldsymbol{v} \, d\boldsymbol{x}\,dt
= \frac{T}{2} \sum_{k=1}^{\infty} k \omega (\sigma 
 \boldsymbol{y}_k,\boldsymbol{v}_k)_{\Omega}
\end{align}
and
\begin{align}
\label{def:weighted-half-time-product:perp}
 \big(\sigma 
 \partial^{1/2}_t \boldsymbol{y}, \partial^{1/2}_t \boldsymbol{v}^{\perp}\big) :=
 \int_{Q} \sigma \partial_t^{1/2} \boldsymbol{y} \cdot
\partial_t^{1/2} \boldsymbol{v}^{\perp} \, d\boldsymbol{x}\,dt
= \frac{T}{2} \sum_{k=1}^{\infty} k \omega (\sigma 
 \boldsymbol{y}_k,\boldsymbol{v}_k^{\perp})_{\Omega}.
\end{align}
$\boldsymbol{H}^{\emph{\textbf{curl}},\frac{1}{2}}_{per}(Q)$-seminorm and norm are defined as
\begin{align*}
 |\boldsymbol{v}|_{\boldsymbol{H}^{\emph{\textbf{curl}},\frac{1}{2}}}^2 :=
 T\|\emph{\textbf{curl}} \, \boldsymbol{v}_0^c\|_{\Omega}^2 + \frac{T}{2} \sum_{k=1}^{\infty} \Big( k\omega \|\boldsymbol{v}_k\|_{\Omega}^2 + \|\emph{\textbf{curl}} \, \boldsymbol{v}_k\|_{\Omega}^2 \Big)
 \end{align*}
 and
 \begin{align*}
 \|\boldsymbol{v}\|_{\boldsymbol{H}^{\emph{\textbf{curl}},\frac{1}{2}}}^2 :=
 T(\|\boldsymbol{v}_0^c\|_{\Omega}^2 + \|\emph{\textbf{curl}} \, \boldsymbol{v}_0^c\|_{\Omega}^2) + \frac{T}{2} \sum_{k=1}^{\infty} \Big( (1+k\omega) \|\boldsymbol{v}_k\|_{\Omega}^2 + \|\emph{\textbf{curl}} \, \boldsymbol{v}_k\|_{\Omega}^2] \Big).
\end{align*}
\end{definition}
The following lemma can be found in \cite{LangerWolfmayr:2013, LangerRepinWolfmayr:2015} for the scalar case. We introduce here the vector-valued version.
\begin{lemma}
 \label{lemma:H11/2IdentiesAndOrthogonalities}
 The identities
 \begin{align}
 \label{equation:H11/2identities}
 \begin{aligned}
  \big(\sigma \partial_t^{1/2} \boldsymbol{y},\partial_t^{1/2} \boldsymbol{v} \big)  =
  \big(\sigma \partial_t \boldsymbol{y},\boldsymbol{v}^{\perp} \big)  \quad \mbox{ and } \quad 
  \big(\sigma \partial_t^{1/2} \boldsymbol{y},\partial_t^{1/2} \boldsymbol{v}^{\perp} \big)  =
  \big(\sigma \partial_t \boldsymbol{y},\boldsymbol{v} \big) 
 \end{aligned}
 \end{align}
 are valid 
 for all $\boldsymbol{y} \in \boldsymbol{H}^{0,1}_{per}(Q)$
 and $\boldsymbol{v} \in \boldsymbol{H}^{0,\frac{1}{2}}_{per}(Q)$.
\end{lemma}
The orthogonality relations
\begin{align}
\label{equation:orthorelation}
\begin{aligned}
 &\big(\sigma \partial_t \boldsymbol{y},\boldsymbol{y}\big)  = 0 \quad \mbox{and} \quad (\sigma \boldsymbol{y}^{\perp},\boldsymbol{y})  = 0 \qquad \forall \,
 \boldsymbol{y} \in \boldsymbol{H}^{0,1}_{per}(Q), \\
 &\big(\sigma \partial^{1/2}_t \boldsymbol{y},\partial^{1/2}_t \boldsymbol{y}^{\perp}\big) = 0
 \quad \mbox{and} \quad
 \big(\nu \, \textbf{curl} \, \boldsymbol{y}, \textbf{curl} \, \boldsymbol{y}^{\perp}\big) = 0
 \quad \forall \, \boldsymbol{y} \in \boldsymbol{H}^{\textbf{curl},\frac{1}{2}}_{per}(Q)
\end{aligned}
\end{align}
hold with
$\big(\nu \, \textbf{curl} \, \boldsymbol{y}, \textbf{curl} \, \boldsymbol{y}^{\perp}\big) := \sum_{k=1}^{\infty} (\nu \, \textbf{curl} \, \boldsymbol{y}_k, \textbf{curl} \, \boldsymbol{y}_k^\perp)_\Omega$
and
$\textbf{curl} \, \boldsymbol{y}_k^\perp
:= (- \textbf{curl} \, \boldsymbol{y}_k^s, \textbf{curl} \, \boldsymbol{y}_k^c)^T$ for all $k \in \mathbb{N}$.
The Friedrichs inequality for $\boldsymbol{H}^{\textbf{curl}}(\Omega)$ (see, e.g., \cite{GiraultRaviart:1986}) holds also for functions represented by their Fourier series as follows
\begin{align}
\label{inequality:Friedrichs:FourierSpace}
 \begin{aligned}
 \|\textbf{curl} \, \boldsymbol{v}\| ^2 &= \int_{Q} |\textbf{curl} \, \boldsymbol{v}|^2 \, d\boldsymbol{x}\,dt
 = T \, \|\textbf{curl} \, \boldsymbol{v}_0^c\|_{\Omega}^2
 + \frac{T}{2} \sum_{k=1}^\infty \|\textbf{curl} \, \boldsymbol{v}_k\|_{\Omega}^2 \\
 &\geq \frac{1}{{C_F^{\text{curl}}}^2} \left(T \, \|\boldsymbol{v}_0^c\|_{\Omega}^2
 + \frac{T}{2} \sum_{k=1}^\infty \|\boldsymbol{v}_k\|_{\Omega}^2 \right) = \frac{1}{{C_F^{\text{curl}}}^2} \|\boldsymbol{v}\| ^2
 \end{aligned}
\end{align}
for weakly divergent functions $\boldsymbol{v} \in \boldsymbol{H}^{\textbf{curl},0}(Q)$, where ${C_F^{\text{curl}}} > 0$ is a constant depending only on the domain $\Omega$.
The weakly divergence-free condition for the adjoint state $\boldsymbol{p}$ can be stated as
\begin{align}
\label{equation:WeakDivCond}
 \int_{Q} \boldsymbol{p} \cdot \nabla v \, d\boldsymbol{x}\,dt = 0 \qquad \forall \, v \in H^{1,0}_0(Q).
\end{align}
We include the gauging condition by introducing the following definitions.
\begin{definition}
We define the spaces
\begin{align*}
\boldsymbol{W} :=
\{ &\boldsymbol{v} \in \boldsymbol{H}^{\textbf{curl}}(\Omega): \exists \, \psi \in H^{1}(\Omega):  \boldsymbol{v} = \nabla \psi, (\psi,1)_{\Omega} = 0, \, \psi |_{\Gamma} = c, \, c \in \mathbb{R} \},  \\
\boldsymbol{H}^{\textbf{curl}|_0,0}(Q) := \{ &\boldsymbol{v} \in \boldsymbol{H}^{\textbf{curl},0}(Q): (\sigma \boldsymbol{y}(t), \boldsymbol{v})_{\Omega} = 0 \,\, \forall \, \boldsymbol{v} \in \boldsymbol{W} \, \text{for a.e.} \, t \in (0,T) \}, \\
  \boldsymbol{H}^{\textbf{curl}|_0,\frac{1}{2}}_{0,per}(Q)
:= \{ &\boldsymbol{v} \in \boldsymbol{H}^{\textbf{curl},\frac{1}{2}}_{0,per}(Q):
(\sigma \boldsymbol{y}(t), \boldsymbol{v})_{\Omega} = 0 
 \,\, \forall \, \boldsymbol{v} \in \boldsymbol{W} \, \text{for a.e.} \, t \in (0,T) \},
\end{align*}
Hence, the latter two are the gauged subspaces of $\boldsymbol{H}^{\textbf{curl},0}(Q)$
and $\boldsymbol{H}^{\textbf{curl},\frac{1}{2}}_{0,per}(Q)$, respectively.
\end{definition}
\begin{remark}
The divergence-free property 
\eqref{equation:divergencefree:yd}, which is passed over to 
$\boldsymbol{u}$ is valid in the frequency domain, so for all modes $k \in \mathbb{N}_0$. This leads to
\begin{align*}
0 = (\boldsymbol{u}(t), \nabla \psi)_{\Omega} = (\boldsymbol{u}_0^c, \nabla \psi)_{\Omega}
+ \sum_{k=0}^{\infty} [(\boldsymbol{u}_k^c, \nabla \psi)_{\Omega} \cos(k \omega t) 
+ (\boldsymbol{u}_k^s, \nabla \psi)_{\Omega} \sin(k \omega t)]
\end{align*}
for all $t \in (0,T)$.
Due to the orthogonality of cosine and sine functions, we immediately obtain that the single Fourier coefficients are weakly divergence-free.
\end{remark}
The weakly divergence-free condition \eqref{equation:divergencefree:yd} for the desired state reads for the Fourier coefficients as follows
\begin{align}
\label{equation:divergencefree:yd:k}
 \int_\Omega \boldsymbol{y_d}_k \cdot \nabla v \, d\boldsymbol{x} = 0 \qquad \forall \, v \in H^1_0(\Omega)
\end{align}
for all $k \in \mathbb{N}$
and
\begin{align}
\label{equation:divergencefree:yd:0}
 \int_\Omega \boldsymbol{y_d}_0^c \cdot \nabla v \, d\boldsymbol{x} = 0 \qquad \forall \, v \in H^1_0(\Omega)
\end{align}
for $k=0$.
\begin{remark}
The gauging condition is essential in the case $k=0$.
The coefficient $\boldsymbol{y}_0^c$ is constant and a non-unique contribution to the solution $\boldsymbol{y}$ for the eddy current problem corresponding to the forward problem 
only. However, this is not the case for the full optimal control problem.
The condition is always redundant for the cases $k \in \mathbb{N}$, since here it is passed from the right-hand side $\boldsymbol{u}$ to $\boldsymbol{y}_k$ for $k \in \mathbb{N}$ as non-constant unique contributions.
\end{remark}

\section{Weak space-time variational formulations}
\label{Sec4:WSTVF}

In order to derive functional type error estimates, we present
variational formulations for the forward problem 
\eqref{equation:forwardpde:pde}-\eqref{equation:forwardpde:div} and the optimality system \eqref{equation:KKTSysClassical:1}-\eqref{equation:KKTSysClassical:6}. 
We start with the forward problem. 
In this case, the function $\boldsymbol{u}$ takes the role of the given data. 
The problem reads as follows:
Given $\boldsymbol{u} \in \boldsymbol{L^2}(Q)$ fulfilling
\eqref{equation:divergencefree:yd}, find $\boldsymbol{y} \in \boldsymbol{H}^{\textbf{curl},\frac{1}{2}}_{0,per}(Q)$ such that
\begin{align}
\label{equation:STVF11/2}
\int_{Q} \left(\sigma \partial_t^{1/2} \boldsymbol{y} \cdot
\partial_t^{1/2} \boldsymbol{v}^{\perp}
  + \nu \, \textbf{curl} \, \boldsymbol{y} \cdot \textbf{curl} \, \boldsymbol{v} \right)\, d\boldsymbol{x}\,dt
	= \int_{Q} \boldsymbol{u}\cdot\boldsymbol{v} \, d\boldsymbol{x}\,dt
\end{align}
for all $\boldsymbol{v} \in \boldsymbol{H}^{\textbf{curl},\frac{1}{2}}_{0,per}(Q)$.
The functions are expanded into Fourier series.
This is a natural approach because of the time-periodicity condition.
Lemma \ref{lemma:STBFinfsupsupsup} yields existence and uniqueness of a solution for variational problem (\ref{equation:STVF11/2}) by applying Babu\v{s}ka-Aziz' theorem (see \cite{Babuska:1971} and \cite{BabuskaAziz:1972}).

\begin{lemma}
\label{lemma:STBFinfsupsupsup}
 The following inf-sup and sup-sup conditions are fulfilled
  \begin{align}
  \label{inequality:STBFinfsupsupsup:Seminorm}
  \underline{c} |\boldsymbol{y}|_{\boldsymbol{H}^{\emph{\textbf{curl}},\frac{1}{2}}} \leq  \sup_{0 \not= \boldsymbol{v} \in \boldsymbol{H}^{\textbf{curl},\frac{1}{2}}_{0,per}(Q)}    \frac{a(\boldsymbol{y},\boldsymbol{v})}{|\boldsymbol{v}|_{\boldsymbol{H}^{\emph{\textbf{curl}},\frac{1}{2}}}} \leq  \overline{c} |\boldsymbol{y}|_{\boldsymbol{H}^{\emph{\textbf{curl}},\frac{1}{2}}}
 \end{align} 
 for all $\boldsymbol{y} \in \boldsymbol{H}^{\textbf{curl},\frac{1}{2}}_{0,per}(Q)$, where $\underline{c}$ and $\overline{c}$ are positive constants depending only on the maximum and minimum values of the conductivity and reluctivity parameters: $\underline{c} = \min\{\underline{\nu},\underline{\sigma}\}/\sqrt{2}$ and $\overline{c} = \max\{\overline{\sigma},\overline{\nu}\}$.
 The space-time bilinear form is given as follows
 \begin{align*}
  a(\boldsymbol{y},\boldsymbol{v}) = \int_{Q} \left(\sigma \partial_t^{1/2} \boldsymbol{y} \cdot \partial_t^{1/2} \boldsymbol{v}^{\perp}  + \nu \, \emph{\textbf{curl}} \, \boldsymbol{y} \cdot \emph{\textbf{curl}} \, \boldsymbol{v} \right)\, d\boldsymbol{x}\,dt.
 \end{align*}
\begin{proof}
Applying triangle and Cauchy-Schwarz inequalities yields the upper estimate
 \begin{align*}
   |a(\boldsymbol{y},\boldsymbol{v})| &=
   \Big|\int_{Q} \left(\sigma \partial_t^{1/2} \boldsymbol{y} \cdot
  \partial_t^{1/2} \boldsymbol{v}^{\perp}
  + \nu \, \textbf{curl} \, \boldsymbol{y} \cdot
  \textbf{curl} \, \boldsymbol{v} \right)\, d\boldsymbol{x}\,dt\Big| \\
   &\leq \overline{\sigma} \Big|\int_{Q} \partial_t^{1/2} \boldsymbol{y} \cdot
  \partial_t^{1/2} \boldsymbol{v}^{\perp} \, d\boldsymbol{x}\,dt \Big|
   + \overline{\nu} \Big|\int_{Q}  \textbf{curl} \, \boldsymbol{y} \cdot
  \textbf{curl} \, \boldsymbol{v} \, d\boldsymbol{x}\,dt\Big| \\
   &\leq \overline{\sigma} \big\|\partial^{1/2}_t \boldsymbol{y}\big\| \big\|\partial^{1/2}_t \boldsymbol{v}\big\|
   + \overline{\nu} \| \textbf{curl} \, \boldsymbol{y}\| \| \textbf{curl} \, \boldsymbol{v}\| \\
   &\leq \max\{\overline{\sigma},\overline{\nu}\} | \boldsymbol{y}|_{\boldsymbol{H}^{\emph{\textbf{curl}},\frac{1}{2}}}
    | \boldsymbol{v}|_{\boldsymbol{H}^{\emph{\textbf{curl}},\frac{1}{2}}}
 \end{align*}
deducing the constant $\overline{c} = \max\{\overline{\sigma},\overline{\nu}\}$.
For the lower estimate, we choose the test function $\boldsymbol{v} = \boldsymbol{y} - \boldsymbol{y}^{\perp}$ and
apply the $\sigma$- and $\nu$-weighted orthogonality relations (\ref{equation:orthorelation}).
We obtain the inequalities
  \begin{align*}
   a(\boldsymbol{y},\boldsymbol{y}) &= 
   \int_{Q} \left(\sigma \partial_t^{1/2} \boldsymbol{y} \cdot
  \partial_t^{1/2} \boldsymbol{y}^{\perp}
  + \nu \, \textbf{curl} \, \boldsymbol{y} \cdot
  \textbf{curl} \, \boldsymbol{y} \right)\, d\boldsymbol{x}\,dt \\
   &= \int_{Q} \nu \, \textbf{curl} \, \boldsymbol{y} \cdot
  \textbf{curl} \, \boldsymbol{y} \, d\boldsymbol{x}\,dt 
   \geq \underline{\nu} \, \big\| \textbf{curl} \, \boldsymbol{y} \big\|^2
  \end{align*}
  and
  \begin{align*}
   a(\boldsymbol{y},-\boldsymbol{y}^{\perp}) &= 
   \int_{Q} \left(\sigma \partial_t^{1/2} \boldsymbol{y} \cdot
  \partial_t^{1/2} \boldsymbol{y}
  - \nu \, \textbf{curl} \, \boldsymbol{y} \cdot
  \textbf{curl} \, \boldsymbol{y}^{\perp} \right)\, d\boldsymbol{x}\,dt \\
   &= \int_{Q} \sigma \partial_t^{1/2} \boldsymbol{y} \cdot
  \partial_t^{1/2} \boldsymbol{y} \, d\boldsymbol{x}\,dt
   \geq \underline{\sigma} \, \big\| \partial_t^{1/2} \boldsymbol{y} \big\|^2,
  \end{align*}
  both leading to 
  \begin{align*}
   a(\boldsymbol{y},\boldsymbol{y}-\boldsymbol{y}^{\perp})
   &\geq \underline{\nu} \, \big\| \textbf{curl} \, \boldsymbol{y} \big\|^2
   + \underline{\sigma} \, \big\| \partial_t^{1/2} \boldsymbol{y} \big\|^2 \\
   &\geq \min\{\underline{\nu},\underline{\sigma}\} \left(\big\| \textbf{curl} \, \boldsymbol{y} \big\|^2
   + \big\| \partial_t^{1/2} \boldsymbol{y} \big\|^2 \right).
  \end{align*}
  Defining $\underline{c} = \min\{\underline{\nu},\underline{\sigma}\}/\sqrt{2}$,
  this provides the inf-sup condition
   \begin{align*}
  \sup_{0 \not= \boldsymbol{v} \in \boldsymbol{H}^{\textbf{curl},\frac{1}{2}}_{0,per}(Q)}
    \frac{a(\boldsymbol{y},\boldsymbol{v})}{|\boldsymbol{v}|_{\boldsymbol{H}^{\emph{\textbf{curl}},\frac{1}{2}}}}  
   \geq \frac{a(\boldsymbol{y},\boldsymbol{y}-\boldsymbol{y}^{\perp})}{|\boldsymbol{y}-\boldsymbol{y}^{\perp}|_{\boldsymbol{H}^{\emph{\textbf{curl}},\frac{1}{2}}}}
   \geq \underline{c} |\boldsymbol{y}|_{\boldsymbol{H}^{\emph{\textbf{curl}},\frac{1}{2}}}.
   \end{align*} 
 \end{proof}
\end{lemma}

The space-time variational formulation of the optimality system \eqref{equation:KKTSysClassical:1}-\eqref{equation:KKTSysClassical:6} is obtained in the same way as for the forward problem.
It is stated as follows:
Given $\boldsymbol{y_d} \in \boldsymbol{L^2}(Q)$ fulfilling
\eqref{equation:divergencefree:yd}, find $\boldsymbol{y}, \boldsymbol{p} \in \boldsymbol{H}^{\emph{\textbf{curl}},\frac{1}{2}}_{0,per}(Q)$ such that
\begin{align}
 \label{equation:KKTSysSTVF:1}
  \int_{Q} \Big( \boldsymbol{y}\cdot\boldsymbol{v} 
  - \nu \, \textbf{curl} \, \boldsymbol{p} \cdot \textbf{curl} \, \boldsymbol{v}
  + \sigma \partial_t^{1/2} &\boldsymbol{p} \cdot
  \partial_t^{1/2} \boldsymbol{v}^{\perp} \Big)\, d\boldsymbol{x}\,dt 
      = \int_{Q} \boldsymbol{y_d}\cdot\boldsymbol{v}\,d\boldsymbol{x}\,dt, \\
 \label{equation:KKTSysSTVF:2}
  \int_{Q} \Big( \nu \, \textbf{curl} \, \boldsymbol{y} \cdot \textbf{curl} \, \boldsymbol{q}
      + \sigma \partial_t^{1/2} \boldsymbol{y} \cdot
  \partial_t^{1/2} &\boldsymbol{q}^{\perp} 
  + \frac{1}{\alpha} \, \boldsymbol{p}\cdot\boldsymbol{q} \Big)\,d\boldsymbol{x}\,dt = 0
\end{align}
for all $\boldsymbol{v}, \boldsymbol{q} \in \boldsymbol{H}^{\emph{\textbf{curl}},\frac{1}{2}}_{0,per}(Q)$.
Similarly, we derive inf-sup and sup-sup conditions for the bilinear form
 \begin{align}
 \label{definition:KKTSysSTVF}
 \begin{aligned}
  \mathcal{B}((\boldsymbol{y},\boldsymbol{p}),(\boldsymbol{v},\boldsymbol{q}))
  = \int_{Q} \Big( &\boldsymbol{y}\cdot\boldsymbol{v} 
  - \nu 
   \, \textbf{curl} \, \boldsymbol{p} \cdot \textbf{curl} \, \boldsymbol{v}
  + \sigma 
  \partial_t^{1/2} \boldsymbol{p} \cdot
  \partial_t^{1/2} \boldsymbol{v}^{\perp} \\
  &+ \nu 
  \, \textbf{curl} \, \boldsymbol{y} \cdot \textbf{curl} \, \boldsymbol{q}
      + \sigma 
      \partial_t^{1/2} \boldsymbol{y} \cdot
  \partial_t^{1/2} \boldsymbol{q}^{\perp} 
  + \frac{1}{\alpha} \, \boldsymbol{p}\cdot\boldsymbol{q} \Big)\,d\boldsymbol{x}\,dt. 
 \end{aligned}
 \end{align}
\begin{lemma}
\label{lemma:STBFinfsupsupsup:KKT}
 The space-time bilinear form \eqref{definition:KKTSysSTVF}
 fulfills the following inf-sup and sup-sup conditions:
  \begin{align}
  \label{inequality:KKTSysSTVFinfsupsupsup}
  \underline{c} \|(\boldsymbol{y},\boldsymbol{p})\|_{\boldsymbol{H}^{\emph{\textbf{curl}},\frac{1}{2}}} \leq
  \sup_{0 \not= (\boldsymbol{v},\boldsymbol{q}) \in (\boldsymbol{H}^{\textbf{curl},\frac{1}{2}}_{0,per}(Q))^2}
    \frac{\mathcal{B}((\boldsymbol{y},\boldsymbol{p}),(\boldsymbol{v},\boldsymbol{q}))}{
    \|(\boldsymbol{v},\boldsymbol{q})\|_{\boldsymbol{H}^{\emph{\textbf{curl}},\frac{1}{2}}}} \leq
  \overline{c} \|(\boldsymbol{y},\boldsymbol{p})\|_{\boldsymbol{H}^{\emph{\textbf{curl}},\frac{1}{2}}}
 \end{align}
 for all $(\boldsymbol{y},\boldsymbol{p}) \in (\boldsymbol{H}^{\textbf{curl},\frac{1}{2}}_{0,per}(Q))^2$, where $\underline{c} = (1+2\max\{\alpha,\frac{1}{\alpha}\})^{-1/2}(\min\{\frac{1}{\sqrt{\alpha}},\underline{\nu},\underline{\sigma}\} \min\{\sqrt{\alpha},\frac{1}{\sqrt{\alpha}}\})$ and $\overline{c} = \max\{1,\frac{1}{\alpha},\overline{\nu},\overline{\sigma}\}$ are positive constants.
\begin{proof}
Applying triangle and Cauchy-Schwarz inequalities yields the upper estimate
\begin{align*}
  \big|\mathcal{B}((\boldsymbol{y},\boldsymbol{p}),(\boldsymbol{v},\boldsymbol{q}))\big|
  = &\, \Big| \int_{Q} \Big( \boldsymbol{y}\cdot\boldsymbol{v} 
  - \nu 
   \, \textbf{curl} \, \boldsymbol{p} \cdot \textbf{curl} \, \boldsymbol{v}
  + \sigma 
  \partial_t^{1/2} \boldsymbol{p} \cdot
  \partial_t^{1/2} \boldsymbol{v}^{\perp} \\
  &\qquad \qquad \, + \nu 
  \, \textbf{curl} \, \boldsymbol{y} \cdot \textbf{curl} \, \boldsymbol{q}
      + \sigma 
      \partial_t^{1/2} \boldsymbol{y} \cdot
  \partial_t^{1/2} \boldsymbol{q}^{\perp} 
  + \frac{1}{\alpha} \, \boldsymbol{p}\cdot\boldsymbol{q} \Big)\,d\boldsymbol{x}\,dt  \Big| \\
  \leq &\, \|\boldsymbol{y}\| \|\boldsymbol{v}\|
      + \overline{\nu} \, \|\textbf{curl} \, \boldsymbol{p}\| \|\textbf{curl} \, \boldsymbol{v}\|
      + \overline{\sigma} \, \big\|\partial^{1/2}_t \boldsymbol{p}\big\| \big\|\partial^{1/2}_t \boldsymbol{v}\big\| \\
      &+ \overline{\nu} \, \|\textbf{curl} \, \boldsymbol{y}\| \|\textbf{curl} \, \boldsymbol{q}\|
      + \overline{\sigma} \, \big\|\partial^{1/2}_t \boldsymbol{y}\big\| \big\|\partial^{1/2}_t \boldsymbol{q}\big\|
      + \frac{1}{\alpha} \, \|\boldsymbol{p}\| \|\boldsymbol{q}\| \\
   \leq &\, \overline{c} \|(\boldsymbol{y},\boldsymbol{p})\|_{\boldsymbol{H}^{\emph{\textbf{curl}},\frac{1}{2}}} \|(\boldsymbol{v},\boldsymbol{q})\|_{\boldsymbol{H}^{\emph{\textbf{curl}},\frac{1}{2}}}
\end{align*}
with $\overline{c} = \max\{1,\frac{1}{\alpha},\overline{\nu},\overline{\sigma}\}$.
The lower estimate is proven by choosing
the test function
\begin{align*}
 (\boldsymbol{v},\boldsymbol{q}) = (\boldsymbol{y} - \frac{1}{\sqrt{\alpha}} \boldsymbol{p}
 - \frac{1}{\sqrt{\alpha}} \boldsymbol{p}^\perp,
          \boldsymbol{p} + \sqrt{\alpha} \boldsymbol{y} - \sqrt{\alpha} \boldsymbol{y}^\perp)
\end{align*}
and applying the $\sigma$- and $\nu$-weighted orthogonality relations (\ref{equation:orthorelation}). We obtain the equations
\begin{align*}
 \mathcal{B}((\boldsymbol{y},\boldsymbol{p}),(\boldsymbol{y},\boldsymbol{p})) 
  =&\, \|\boldsymbol{y}\|^2 + \frac{1}{\alpha} \|\boldsymbol{p}\|^2, \\
 \mathcal{B}((\boldsymbol{y},\boldsymbol{p}),(-\frac{1}{\sqrt{\alpha}} \boldsymbol{p}, \sqrt{\alpha} \boldsymbol{y}))
  = &\, \frac{1}{\sqrt{\alpha}} (\nu \, \textbf{curl} \, \boldsymbol{p},\textbf{curl} \, \boldsymbol{p})
  + \sqrt{\alpha} (\nu \, \textbf{curl} \, \boldsymbol{y},\textbf{curl} \, \boldsymbol{y}), \\
 \mathcal{B}((\boldsymbol{y},\boldsymbol{p}),(-\frac{1}{\sqrt{\alpha}} \boldsymbol{p}^\perp,-\sqrt{\alpha} \boldsymbol{y}^\perp)) 
  = &\, \frac{1}{\sqrt{\alpha}} (\sigma \partial_t^{1/2} \boldsymbol{p}, \partial_t^{1/2} \boldsymbol{p})
    + \sqrt{\alpha} (\sigma \partial_t^{1/2} \boldsymbol{y}, \partial_t^{1/2} \boldsymbol{y}),
\end{align*}
leading to the lower estimate
\begin{align*}
 \mathcal{B}((\boldsymbol{y},\boldsymbol{p})&,(\boldsymbol{v},\boldsymbol{q}))
 \geq 
 \min\{\frac{1}{\sqrt{\alpha}},\underline{\nu},\underline{\sigma}\}
 \min\{\sqrt{\alpha},\frac{1}{\sqrt{\alpha}}\}
\|(\boldsymbol{y},\boldsymbol{p})\|_{\boldsymbol{H}^{\emph{\textbf{curl}},\frac{1}{2}}}^2.
\end{align*}
Together with
\begin{align*}
\|(\boldsymbol{v},\boldsymbol{q})\|_{\boldsymbol{H}^{\emph{\textbf{curl}},\frac{1}{2}}} 
 \leq \left(1+2\max\{\alpha,\frac{1}{\alpha}\}\right)^{1/2} \|(\boldsymbol{y},\boldsymbol{p})\|_{\boldsymbol{H}^{\emph{\textbf{curl}},\frac{1}{2}}},
\end{align*}
we derive the inf-sup estimate of \eqref{inequality:KKTSysSTVFinfsupsupsup} as follows
\begin{align*}
  \sup_{0 \not= (\boldsymbol{v},\boldsymbol{q}) \in (\boldsymbol{H}^{\emph{\textbf{curl}},\frac{1}{2}}_{0,per}(Q))^2}
    \frac{\mathcal{B}((\boldsymbol{y},\boldsymbol{p}),(\boldsymbol{v},\boldsymbol{q}))}{\|(\boldsymbol{v},\boldsymbol{q})\|_{\boldsymbol{H}^{\textbf{curl},\frac{1}{2}}}}
   &\geq \frac{\min\{\frac{1}{\sqrt{\alpha}},\underline{\nu},\underline{\sigma}\}
 \min\{\sqrt{\alpha},\frac{1}{\sqrt{\alpha}}\} \|(\boldsymbol{y},\boldsymbol{p})\|_{\boldsymbol{H}^{\textbf{curl},\frac{1}{2}}}^2}{
 \sqrt{1+2\max\{\alpha,\frac{1}{\alpha}\}} \|(\boldsymbol{y},\boldsymbol{p})\|_{\boldsymbol{H}^{\textbf{curl},\frac{1}{2}}}}
\end{align*}
and defining $\underline{c} = (1+2\max\{\alpha,\frac{1}{\alpha}\})^{-1/2}(\min\{\frac{1}{\sqrt{\alpha}},\underline{\nu},\underline{\sigma}\}\min\{\sqrt{\alpha},\frac{1}{\sqrt{\alpha}}\})$.
\end{proof}
\end{lemma}

\section{Discretization}
\label{Sec5:MFEM}

We discretize the problem by the multiharmonic finite element method.
All functions in the space-time variational problems are expanded into Fourier series and the problem setting is shifted to the frequency domain.
Using the linearity of the problem and $L^2$-orthogonality relations of cosine and sine functions, 
we derive variational problems for all the Fourier modes $k \in \mathbb{N}$.
For the optimality system:
Given $\boldsymbol{y_d}_k \in (\boldsymbol{L^2}(\Omega))^2$ satisfying \eqref{equation:divergencefree:yd:k}, 
find $\boldsymbol{y}_k = (\boldsymbol{y}_{k}^c, \boldsymbol{y}_{k}^s)^T, \boldsymbol{p}_k = (\boldsymbol{p}_{k}^c, \boldsymbol{p}_{k}^s)^T \in \mathbb{V} = V \times V = (\boldsymbol{H}^{\textbf{curl}}_0(\Omega))^2$ such that
\begin{align}
 \label{equation:MultiAnsVFBlock:1}
 &\int_{\Omega} \big( \boldsymbol{y}_k \cdot \boldsymbol{v}_k
  - \nu \, 
  \textbf{curl} \, \boldsymbol{p}_k \cdot \textbf{curl} \, \boldsymbol{v}_k
  + k \omega \sigma 
  \, \boldsymbol{p}_k \cdot \boldsymbol{v}_k^{\perp} \big)\,d\boldsymbol{x}
  = \int_{\Omega} \boldsymbol{y_d}_k \cdot \boldsymbol{v}_k\,d\boldsymbol{x}, \\
 \label{equation:MultiAnsVFBlock:2}
 &\int_{\Omega} \big( \nu 
 \, \textbf{curl} \, \boldsymbol{y}_k \cdot \textbf{curl} \, \boldsymbol{q}_k
  + k \omega \sigma 
  \,\boldsymbol{y}_k \cdot \boldsymbol{q}_k^{\perp}
  + \frac{1}{\alpha} \, \boldsymbol{p}_k \cdot \boldsymbol{q}_k \big)\,d\boldsymbol{x} = 0
\end{align}
for all $\boldsymbol{v}_k, \boldsymbol{q}_k \in \mathbb{V}$.
The variational problem for the case of $k = 0$ is given by:
Given $\boldsymbol{y_d}_0^c \in \boldsymbol{L^2}(\Omega)$ satisfying \eqref{equation:divergencefree:yd:0},
find $\boldsymbol{y}_0^c, \boldsymbol{p}_0^c \in V = \boldsymbol{H}^{\textbf{curl}}_0(\Omega)$ such that
\begin{align}
 \label{equation:MultiAnsVFBlock0:1}
 &\int_{\Omega} \big( \boldsymbol{y}_0^c \cdot \boldsymbol{v}_0^c
 - \nu \, 
  \textbf{curl} \, \boldsymbol{p}_0^c \cdot \textbf{curl} \, \boldsymbol{v}_0^c \big)\,d\boldsymbol{x}
 = \int_{\Omega}  \boldsymbol{y_d}^c_0 \cdot \boldsymbol{v}_0^c\,d\boldsymbol{x} \\
 \label{equation:MultiAnsVFBlock0:2}
 &\int_{\Omega} \big( \nu \, 
  \textbf{curl} \, \boldsymbol{y}_0^c \cdot \textbf{curl} \, \boldsymbol{q}_0^c + \frac{1}{\alpha} 
  \, \boldsymbol{p}_0^c \cdot \boldsymbol{q}_0^c \big)\,d\boldsymbol{x} = 0
\end{align}
for all $\boldsymbol{v}_0^c, \boldsymbol{q}_0^c \in V$.
Next the Fourier series are truncated. We denote by $N$ the truncation index.
We use the finite element functions 
$\boldsymbol{y}_{kh} = (\boldsymbol{y}_{kh}^c, \boldsymbol{y}_{kh}^s)^T$,
$\boldsymbol{p}_{kh} = (\boldsymbol{p}_{kh}^c, \boldsymbol{p}_{kh}^s)^T
\in \mathbb{V}_h = V_h \times V_h \subset \mathbb{V}$
to approximate the Fourier coefficients
$\boldsymbol{y}_k = (\boldsymbol{y}_k^c, \boldsymbol{y}_k^s)^T$,
$\boldsymbol{p}_k = (\boldsymbol{p}_k^c, \boldsymbol{p}_k^s)^T \in \mathbb{V}$
with the finite element spaces $\mathbb{V}_h = V_h \times V_h$ and
$V_h = \mbox{span} \{\boldsymbol{\phi}_1, ..., \boldsymbol{\phi}_{N_h}\}$,
the discretization parameter $h$ and the dimension of $V_h$ given by $N_h = \mathcal{O}(h^{-3})$. 
We use the N\'{e}d\'{e}lec (edge) basis functions of lowest order 
for approximating the space $V = \boldsymbol{H}^{\textbf{curl}}_0(\Omega)$
(see \cite{Nedelec:1980} and \cite{Nedelec:1986}).
We derive the system of linear equations
\begin{align}
 \label{equation:MultiFESysBlock}
 \left( \begin{array}{cccc}
     \boldsymbol{M_h}  &  0 & -\boldsymbol{K_h} & k \omega \boldsymbol{M_{h,\sigma}} \\
     0  &  \boldsymbol{M_h} & -k \omega \boldsymbol{M_{h,\sigma}} & -\boldsymbol{K_h} \\
     -\boldsymbol{K_h}  &  -k \omega \boldsymbol{M_{h,\sigma}} & -\alpha^{-1} \boldsymbol{M_h} & 0 \\
     k \omega \boldsymbol{M_{h,\sigma}}  &  -\boldsymbol{K_h} & 0 & -\alpha^{-1} \boldsymbol{M_h} \end{array} \right) \left( \begin{array}{c}
     \underline{\boldsymbol{y}}_k^c \\
     \underline{\boldsymbol{y}}_k^s \\
     \underline{\boldsymbol{p}}_k^c \\
     \underline{\boldsymbol{p}}_k^s \end{array} \right) = \left( \begin{array}{c}
     {\underline{\boldsymbol{y}}_{\boldsymbol{d}}^c}_k \\
     {\underline{\boldsymbol{y}}_{\boldsymbol{d}}^s}_k \\
     0 \\
     0 \end{array} \right),
\end{align}
corresponding to \eqref{equation:MultiAnsVFBlock:1}-\eqref{equation:MultiAnsVFBlock:2}.
(Weighted) mass matrices and stiffness matrix are given by
\begin{align*}
\begin{aligned}
 (\boldsymbol{M_h})_{ij} =(\boldsymbol{\varphi}_i, \boldsymbol{\varphi}_j)_{\Omega}, \hspace{0.3cm}
 (\boldsymbol{M_{h,\sigma}})_{ij} &= (\sigma \boldsymbol{\varphi}_i, \boldsymbol{\varphi}_j)_{\Omega}, \hspace{0.3cm}
 (\boldsymbol{K_h})_{ij} &= (\nu \, \textbf{curl} \, \boldsymbol{\varphi}_i, \textbf{curl} \, \boldsymbol{\varphi}_j)_{\Omega},
\end{aligned}
\end{align*}
where $i,j = 1,...,N_h$.
The system of linear equations for the Fourier mode $k=0$ is given by
\begin{align}
 \label{equation:MultiFESysBlock0}
 \left( \begin{array}{cc}
     \boldsymbol{M_h}  &  -\boldsymbol{K_h} \\
     -\boldsymbol{K_h}  &  - \alpha^{-1} \boldsymbol{M_h} \end{array} \right) \left( \begin{array}{c}
     \underline{\boldsymbol{y}}_0^c \\
     \underline{\boldsymbol{p}}_0^c \end{array} \right) = \left( \begin{array}{c}
     {\underline{\boldsymbol{y}}_d^c}_0 \\
     0 \end{array} \right).
\end{align}
Adding up the solutions of linear systems (\ref{equation:MultiFESysBlock}) and (\ref{equation:MultiFESysBlock0}), which provide approximations for the Fourier coefficients of state and adjoint state, and inserting them in the truncated Fourier series yields the multiharmonic finite element approximations for the state and adjoint state which are 
\begin{align}
 \label{definition:MultiharmonicFEAnsatzState}
   \boldsymbol{y}_{N h}(\boldsymbol{x},t) &= \sum_{k=0}^N [\boldsymbol{y}_{kh}^c(\boldsymbol{x}) \cos(k \omega t)
      + \boldsymbol{y}_{kh}^s(\boldsymbol{x}) \sin(k \omega t)], \\
 \label{definition:MultiharmonicFEAnsatzStateCostate}
   \boldsymbol{p}_{N h}(\boldsymbol{x},t) &= \sum_{k=0}^N [\boldsymbol{p}_{kh}^c(\boldsymbol{x}) \cos(k \omega t)
      + \boldsymbol{p}_{kh}^s(\boldsymbol{x}) \sin(k \omega t)].
\end{align}
Fast and robust solvers for the saddle point systems (\ref{equation:MultiFESysBlock}) and
(\ref{equation:MultiFESysBlock0}) can be found, e.g., in
\cite{KolmbauerLanger:2012, Kolmbauer:2012:thesis} or in
\cite{AxelssonLukas:2019}.

Similarly we derive the multiharmonic variational problems for the forward problem:
Given $\boldsymbol{u}_k \in (\boldsymbol{L^2}(\Omega))^2$, find $\boldsymbol{y}_k \in \mathbb{V} = V \times V = (\boldsymbol{H}^{\textbf{curl}}_0(\Omega))^2$ such that
\begin{align}
 \label{equation:MultiAnsVFBlock:forward}
 \int_{\Omega} \big( k \omega \sigma \, \boldsymbol{y}_k \cdot \boldsymbol{v}_k^{\perp}
  + \nu \,  \textbf{curl} \, \boldsymbol{y}_k \cdot \textbf{curl} \, \boldsymbol{v}_k \big)\,d\boldsymbol{x}
  = \int_{\Omega} \boldsymbol{u}_k \cdot \boldsymbol{v}_k\,d\boldsymbol{x},
\end{align}
for all $\boldsymbol{v}_k \in \mathbb{V}$.
For the case $k = 0$, we have:
Given $\boldsymbol{u}_0^c \in \boldsymbol{L^2}(\Omega)$, find $\boldsymbol{y}_0^c \in V = \boldsymbol{H}^{\textbf{curl}}_0(\Omega)$ such that
\begin{align}
 \label{equation:MultiAnsVFBlock0:forward}
 &\int_{\Omega} \big( \nu \,  \textbf{curl} \, \boldsymbol{y}_0^c \cdot \textbf{curl} \, \boldsymbol{v}_0^c \big)\,d\boldsymbol{x}
 = \int_{\Omega}  \boldsymbol{u}^c_0 \cdot \boldsymbol{v}_0^c\,d\boldsymbol{x}
\end{align}
for all $\boldsymbol{v}_0^c \in V$.
The multiharmonic finite element discretization for the forward problem leads to the systems of linear equations as follows
\begin{align}
 \label{equation:MultiFESysBlock:forward}
 \left( \begin{array}{cc}
  \boldsymbol{K_h} & k \omega \boldsymbol{M_{h,\sigma}} \\
  -k \omega \boldsymbol{M_{h,\sigma}} & \boldsymbol{K_h} \\ \end{array} \right) 
  \left( \begin{array}{c}
     \underline{\boldsymbol{y}}_k^c \\
     \underline{\boldsymbol{y}}_k^s \end{array} \right) = \left( \begin{array}{c}
     \underline{\boldsymbol{u}}_k^c \\
     \underline{\boldsymbol{u}}_k^s \end{array} \right),
\end{align}
and $\boldsymbol{K_h} \underline{\boldsymbol{y}}_0^c = \underline{\boldsymbol{u}}_0^c$ for $k=0$.

In the next sections, we will present functional a posteriori estimates for two different problems:
first, for the PDE constraint (the time-periodic eddy current (forward) problem), and second,
for the optimal control problem's optimality system.

\section{A posteriori error estimates for the forward problem}
\label{Sec6:FunctionalAPosterioriEstimates:Forwardpde}

We denote by the function $\boldsymbol{\eta}$ an approximation for the state function $\boldsymbol{y}$.
Let $\boldsymbol{\eta} \in \boldsymbol{H}^{\textbf{curl},1}_{0,per}(Q)$.
Note that $\boldsymbol{\eta}$ is arbitrary for now but we will choose
$\boldsymbol{y}_{Nh}$ as $\boldsymbol{\eta}$. 
Note that for the foward problem alone, we have assumed that $\boldsymbol{u} \in \boldsymbol{L^2}(Q)$ is given and fulfills \eqref{equation:divergencefree:yd}.
We want to derive an estimate from above for the error
$\boldsymbol{y} - \boldsymbol{\eta}$ in $\boldsymbol{H}^{\textbf{curl},\frac{1}{2}}_{0,per}(Q)$.
The bilinear form
$a(\boldsymbol{\boldsymbol{y} - \boldsymbol{\eta}},\boldsymbol{v})$ equals
\begin{align}
 \label{problem:STVFAPost:Error}
  \begin{aligned}
 \int_{Q} &\Big( \sigma
 \partial_t^{1/2} (\boldsymbol{y} - \boldsymbol{\eta}) \cdot \partial_t^{1/2} \boldsymbol{v}^{\perp}
  + \nu \, \textbf{curl} \, (\boldsymbol{y} - \boldsymbol{\eta}) \cdot
  \, \textbf{curl} \, \boldsymbol{v} \Big) d\boldsymbol{x}\,dt \\
 &= \int_{Q} \Big(\boldsymbol{u} \cdot \boldsymbol{v} - \sigma
 \partial_t^{1/2} \boldsymbol{\eta} \cdot \partial_t^{1/2} \boldsymbol{v}^{\perp}
  - \nu  \, \textbf{curl} \, \boldsymbol{\eta} \cdot \, \textbf{curl} \, \boldsymbol{v} \Big) \, d\boldsymbol{x}\,dt
 \end{aligned}
\end{align}
for all $\boldsymbol{v} \in \boldsymbol{H}^{\textbf{curl},\frac{1}{2}}_{0,per}(Q)$.
We define the right-hand side of \eqref{problem:STVFAPost:Error} as a linear functional
for $\boldsymbol{v}$
as follows
\begin{align*}
 \mathcal{F}_{\boldsymbol{\eta}}(\boldsymbol{v}) = \int_{Q} \Big(\boldsymbol{u} \cdot \boldsymbol{v} - \sigma
 \partial_t^{1/2} \boldsymbol{\eta} \cdot \partial_t^{1/2} \boldsymbol{v}^{\perp}
  - \nu \, \textbf{curl} \, \boldsymbol{\eta} \cdot \, \textbf{curl} \, \boldsymbol{v} \Big) \, d\boldsymbol{x}\,dt.
\end{align*}
\begin{theorem}
\label{theorem:aposteriorEstimateH11/2Seminorm}
 Let $\boldsymbol{\eta} \in \boldsymbol{H}^{\textbf{curl},1}_{0,per}(Q)$.
 Let the bilinear form $a(\cdot,\cdot)$ fulfill the inf-sup 
 condition in (\ref{inequality:STBFinfsupsupsup:Seminorm}). The error between the exact state $\boldsymbol{y}$ and $\boldsymbol{\eta}$ can be estimated from above as follows
 \begin{align}
  \label{inequality:aposteriorEstimateH11/2Seminorm}
  |\boldsymbol{y}-\boldsymbol{\eta}|_{\boldsymbol{H}^{\textbf{curl},\frac{1}{2}}} \leq \frac{1}{\underline{c}} \left({C_F^{\text{curl}}} \, \|\mathcal{R}_1(\boldsymbol{\eta},\boldsymbol{\tau})\| + \|\mathcal{R}_2(\boldsymbol{\eta},\boldsymbol{\tau})\| \right)
 \end{align}
 with $\boldsymbol{\tau} \in \boldsymbol{H}^{\textbf{curl},0}(Q)$,
 the Friedrichs constant ${C_F^{\text{curl}}}$ as defined in \eqref{inequality:Friedrichs:FourierSpace},
 and the constant $\underline{c} = \frac{1}{\sqrt{2}}\min\{\underline{\nu},\underline{\sigma}\}$.
 The residual functions $\mathcal{R}_1$ and $\mathcal{R}_2$ are defined as follows
 \begin{align}
\label{def:R1R2}
 \mathcal{R}_1( \boldsymbol{\eta},\boldsymbol{\tau})
   = \boldsymbol{u} - \sigma \partial_t \boldsymbol{\eta} - \textbf{curl} \, \boldsymbol{\tau} 
  \qquad \text{ and } \qquad
 \mathcal{R}_2( \boldsymbol{\eta},\boldsymbol{\tau})
   = \boldsymbol{\tau} - \nu  \, \textbf{curl} \, \boldsymbol{\eta}.
\end{align}
 \begin{proof}
 For all $\boldsymbol{\eta} \in \boldsymbol{H}^{\textbf{curl},1}_{0,per}(Q)$ and $\boldsymbol{v} \in \boldsymbol{H}^{\textbf{curl},\frac{1}{2}}_{0,per}(Q)$, the identity
\begin{align}
\label{equation:identityEtaH11/2}
 \big(\sigma \partial_t^{1/2} \boldsymbol{\eta} ,\partial_t^{1/2} \boldsymbol{v}^{\perp} \big)  =
 \big(\sigma \partial_t \boldsymbol{\eta} ,\boldsymbol{v} \big)
\end{align} 
is valid using \eqref{equation:H11/2identities}. Together with
\begin{align}
\label{identity:curl}
 \int_{Q} \textbf{curl} \, \boldsymbol{\tau} \cdot \boldsymbol{v} \, d\boldsymbol{x} \, dt
 = \int_{Q}  \boldsymbol{\tau} \cdot \textbf{curl} \, \boldsymbol{v} \, d\boldsymbol{x} \, dt \qquad \forall \, \boldsymbol{v} \in \boldsymbol{H}^{\textbf{curl},0}_0(Q) \quad \forall \,
\boldsymbol{\tau} \in \boldsymbol{H}^{\textbf{curl},0}(Q),
\end{align}
this yields the functional $\mathcal{F}_{\boldsymbol{\eta}}(\boldsymbol{v})$ to be represented as
\begin{align*}
\mathcal{F}_{\boldsymbol{\eta}}(\boldsymbol{v}) = \int_{Q} \Big((\boldsymbol{u}  - \sigma \partial_t \boldsymbol{\eta} - \textbf{curl} \, \boldsymbol{\tau}) \cdot \boldsymbol{v} + (\boldsymbol{\tau} - \nu  \, \textbf{curl} \, \boldsymbol{\eta}) \cdot \textbf{curl} \, \boldsymbol{v} \Big) \, d\boldsymbol{x} \, dt \qquad \forall \, \boldsymbol{\tau} \in \boldsymbol{H}^{\textbf{curl},0}(Q).
\end{align*}
We obtain the estimate
\begin{align}
\label{inequality:supRHS:CS}
\mathcal{F}_{\boldsymbol{\eta}}(\boldsymbol{v})
 \leq \|\mathcal{R}_1(\boldsymbol{\eta},\boldsymbol{\tau})\|  \|\boldsymbol{v}\|  + \|\mathcal{R}_2(\boldsymbol{\eta},\boldsymbol{\tau})\|  \|\textbf{curl} \, \boldsymbol{v}\|
\end{align}
by applying the Cauchy-Schwarz inequality with the residual functions $\mathcal{R}_1$ and $\mathcal{R}_2$ defined in \eqref{def:R1R2}.
Applying the Friedrichs inequality  for $\boldsymbol{H}^{\textbf{curl}}(\Omega)$ in \eqref{inequality:Friedrichs:FourierSpace} yields the estimate
\begin{align*}
 \mathcal{F}_{\boldsymbol{\eta}}(\boldsymbol{v}) \leq \|\mathcal{R}_1(\boldsymbol{\eta},\boldsymbol{\tau})\|  \|\boldsymbol{v}\| + \|\mathcal{R}_2(\boldsymbol{\eta},\boldsymbol{\tau})\|  \|\textbf{curl} \, \boldsymbol{v}\| \leq 
 \left({C_F^{\text{curl}}} \, \|\mathcal{R}_1(\boldsymbol{\eta},\boldsymbol{\tau})\| + \|\mathcal{R}_2(\boldsymbol{\eta},\boldsymbol{\tau})\| \right) \|\textbf{curl} \, \boldsymbol{v}\|.
\end{align*}
This leads to
\begin{align}
\label{inequality:aposteriorEstimateH11/2Seminorm0Step}
\begin{aligned}
  \sup_{0 \not= \boldsymbol{v} \in \boldsymbol{H}^{\textbf{curl},\frac{1}{2}}_{0,per}(Q)}
    \frac{\mathcal{F}_{\boldsymbol{\eta}}(\boldsymbol{v})}{|\boldsymbol{v}|_{\boldsymbol{H}^{\textbf{curl},\frac{1}{2}}}}
  &\leq {C_F^{\text{curl}}} \, \|\mathcal{R}_1(\boldsymbol{\eta},\boldsymbol{\tau})\| 
 + \|\mathcal{R}_2(\boldsymbol{\eta},\boldsymbol{\tau})\|.
\end{aligned}
\end{align}
Due to \eqref{problem:STVFAPost:Error}, we have
\begin{align*}
\sup_{0 \not= \boldsymbol{v} \in \boldsymbol{H}^{\textbf{curl},\frac{1}{2}}_{0,per}(Q)}    \frac{a(\boldsymbol{y}-\boldsymbol{\eta},\boldsymbol{v})}{|\boldsymbol{v}|_{\boldsymbol{H}^{\emph{\textbf{curl}},\frac{1}{2}}}} = \sup_{0 \not= \boldsymbol{v} \in \boldsymbol{H}^{\textbf{curl},\frac{1}{2}}_{0,per}(Q)}    \frac{\mathcal{F}_{\boldsymbol{\eta}}(\boldsymbol{v})}{|\boldsymbol{v}|_{\boldsymbol{H}^{\textbf{curl},\frac{1}{2}}}}.
\end{align*}
Together with applying the inf-sup condition \eqref{inequality:STBFinfsupsupsup:Seminorm}, we get
  \begin{align*}
  |\boldsymbol{y}-\boldsymbol{\eta}|_{\boldsymbol{H}^{\textbf{curl},\frac{1}{2}}} \leq \frac{1}{\underline{c}}  \sup_{0 \not= \boldsymbol{v} \in \boldsymbol{H}^{\textbf{curl},\frac{1}{2}}_{0,per}(Q)}    \frac{a(\boldsymbol{y}-\boldsymbol{\eta},\boldsymbol{v})}{|\boldsymbol{v}|_{\boldsymbol{H}^{\textbf{curl},\frac{1}{2}}}}  = \frac{1}{\underline{c}}  \sup_{0 \not= \boldsymbol{v} \in \boldsymbol{H}^{\textbf{curl},\frac{1}{2}}_{0,per}(Q)}    \frac{\mathcal{F}_{\boldsymbol{\eta}}(\boldsymbol{v})}{|\boldsymbol{v}|_{\boldsymbol{H}^{\textbf{curl},\frac{1}{2}}}},
 \end{align*}
 which proves \eqref{inequality:aposteriorEstimateH11/2Seminorm}.
 \end{proof}
\end{theorem}
Note that in case of $\mathcal{R}_1(\boldsymbol{\eta},\boldsymbol{\tau}) = 0$ and $\mathcal{R}_2(\boldsymbol{\eta},\boldsymbol{\tau}) = 0$, we have the eddy-current problem \eqref{equation:forwardpde:pde} represented as $\sigma \partial_t \boldsymbol{\eta} + \textbf{curl} \, \boldsymbol{\tau} = \boldsymbol{u}$ and $\boldsymbol{\tau} = \nu \, \textbf{curl} \, \boldsymbol{\eta}$, where $\boldsymbol{\eta} \in \boldsymbol{H}^{\textbf{curl},1}_{0,per}(Q)$ satisfies the boundary and periodicity conditions \eqref{equation:forwardpde:boundary} and \eqref{equation:forwardpde:periodic} providing the solution of problem \eqref{equation:forwardpde:pde}. Hence, $\mathcal{R}_1(\boldsymbol{y},\nu \, \textbf{curl} \, \boldsymbol{y}) = 0$ and $\mathcal{R}_2(\boldsymbol{y},\nu \, \textbf{curl} \, \boldsymbol{y}) = 0$.

We define the upper bound as majorant function
 \begin{align}
  \label{def:PDEconstraint:majorant:seminorm}
  \mathcal{M}_{|\cdot|}^{\oplus}(\boldsymbol{\eta},\boldsymbol{\tau}) = \frac{1}{\underline{c}} \left({C_F^{\text{curl}}} \, \|\mathcal{R}_1(\boldsymbol{\eta},\boldsymbol{\tau})\| + \|\mathcal{R}_2(\boldsymbol{\eta},\boldsymbol{\tau})\| \right)
 \end{align}
which is guaranteed and computable.
The majorant can be estimated by its quadratic representative applying Young's inequality:
 \begin{align}
 \label{maj:quadratic:forward}
  \mathcal{M}^\oplus_{|\cdot|}(\boldsymbol{\eta},\boldsymbol{\tau})^2 \leq 
  \mathcal{M}^\oplus_{|\cdot|}(\beta; \boldsymbol{\eta},\boldsymbol{\tau})^2 
  = \frac{1}{\underline{c}^2} \big({C_F^{\text{curl}}}^2 (1+\beta) \|\mathcal{R}_1(\boldsymbol{\eta},\boldsymbol{\zeta},\boldsymbol{\tau})\|^2 
  + \frac{(1+\beta)}{\beta} \|\mathcal{R}_2(\boldsymbol{\eta},\boldsymbol{\tau})\|^2 \big)
 \end{align}
with the constant $\beta > 0$.
Due to the Friedrichs inequality the seminorm and norm in $\boldsymbol{H}^{\textbf{curl},\frac{1}{2}}$ are equivalent providing the inf-sup and sup-sup conditions
 \begin{align}
  \label{inequality:STBFinfsupsupsup:Norm}
  \underline{c} \|\boldsymbol{y}\|_{\boldsymbol{H}^{\emph{\textbf{curl}},\frac{1}{2}}} \leq  \sup_{0 \not= \boldsymbol{v} \in \boldsymbol{H}^{\textbf{curl},\frac{1}{2}}_{0,per}(Q)}    \frac{a(\boldsymbol{y},\boldsymbol{v})}{\|\boldsymbol{v}\|_{\boldsymbol{H}^{\emph{\textbf{curl}},\frac{1}{2}}}} \leq  \overline{c} \|\boldsymbol{y}\|_{\boldsymbol{H}^{\emph{\textbf{curl}},\frac{1}{2}}}
 \end{align} 
 for all $\boldsymbol{y} \in \boldsymbol{H}^{\textbf{curl},\frac{1}{2}}_{0,per}(Q)$ with the constants $\underline{c} = \min\{\underline{\nu}/(1+{C_F^{\text{curl}}}^2),\underline{\sigma}\}/\sqrt{2}$ and $\overline{c} = \max\{\overline{\sigma},\overline{\nu}\}$, where we use the Friedrichs inequality \eqref{inequality:Friedrichs:FourierSpace} for the lower bound of
\begin{align*}
   a(\boldsymbol{y},\boldsymbol{y}) &= 
   \int_{Q} \left(\sigma(\boldsymbol{x}) \partial_t^{1/2} \boldsymbol{y} \cdot
  \partial_t^{1/2} \boldsymbol{y}^{\perp}
  + \nu(\boldsymbol{x}) \, \textbf{curl} \, \boldsymbol{y} \cdot
  \textbf{curl} \, \boldsymbol{y} \right)\, d\boldsymbol{x}\,dt \\
   &= \int_{Q} \nu(\boldsymbol{x}) \, \textbf{curl} \, \boldsymbol{y} \cdot
  \textbf{curl} \, \boldsymbol{y} \, d\boldsymbol{x}\,dt 
   \geq \underline{\nu} \, \big\| \textbf{curl} \, \boldsymbol{y} \big\|^2
   \geq \frac{\underline{\nu}}{{1+C_F^{\text{curl}}}^2} \big\| \boldsymbol{y} \big\|_{\boldsymbol{H}^{\emph{\textbf{curl}},0}}^2
\end{align*}
following otherwise the proof for \eqref{lemma:STBFinfsupsupsup}.
This provides the majorant 
 \begin{align}
  \label{def:PDEconstraint:majorant:Norm}
  \mathcal{M}_{\|\cdot\|}^{\oplus}(\boldsymbol{\eta},\boldsymbol{\tau}) = \frac{1}{\underline{c}} \sqrt{\|\mathcal{R}_1(\boldsymbol{\eta},\boldsymbol{\tau})\|^2   + \|\mathcal{R}_2(\boldsymbol{\eta},\boldsymbol{\tau})\|^2},
 \end{align}
where we apply
\begin{align*}
\mathcal{F}_{\boldsymbol{\eta}}(\boldsymbol{v})
 \leq \sqrt{\|\mathcal{R}_1(\boldsymbol{\eta},\boldsymbol{\tau})\|^2   + \|\mathcal{R}_2(\boldsymbol{\eta},\boldsymbol{\tau})\|^2}\sqrt{\|\boldsymbol{v}\|^2 + \|\textbf{curl} \, \boldsymbol{v}\|^2}
\end{align*}
in the estimate \eqref{inequality:supRHS:CS}.
The majorants are nonnegative. Only for $\boldsymbol{\eta} = \boldsymbol{y}$ and $\boldsymbol{\tau} = \nu \, \textbf{curl} \, \boldsymbol{y}$, the majorants are zero.

The multiharmonic representation of the majorant \eqref{def:PDEconstraint:majorant:Norm} is given by
\begin{align*}
  \mathcal{M}^\oplus_{\|\cdot\|}(\boldsymbol{\eta},\boldsymbol{\tau}) 
 = &\, \frac{1}{\underline{c}}
  \Big(T \big(\|{\mathcal{R}_1}^c_0(\boldsymbol{\eta}_0^c,\boldsymbol{\tau}_0^c)\|_{\Omega}^2 
  + \|{\mathcal{R}_2}^c_0(\boldsymbol{\eta}_0^c,\boldsymbol{\tau}_0^c)\|_{\Omega}^2 \big) \\
  &+ \frac{T}{2} \sum_{k=1}^N
  \big(\|{\mathcal{R}_1}_k(\boldsymbol{\eta}_k,\boldsymbol{\tau}_k)\|_{\Omega}^2 
  + \|{\mathcal{R}_2}_k(\boldsymbol{\eta}_k,\boldsymbol{\tau}_k)\|_{\Omega}^2 \big)
  + \mathcal{P}_N
  \Big)^{1/2}.
 \end{align*}
The term $\mathcal{P}_N$ is called the remainder term and is fully computable
\begin{align*}
\mathcal{P}_N =
 \frac{T}{2} \sum_{k=N+1}^\infty \|\boldsymbol{u}_k\|_{\Omega}^2
 = \frac{T}{2} \sum_{k=N+1}^\infty \left(\|\boldsymbol{u}_k^c\|_{\Omega}^2
 + \|\boldsymbol{u}_k^s\|_{\Omega}^2\right)
 = \|\boldsymbol{u}- \boldsymbol{u}_N\|
\end{align*}
where $\boldsymbol{u} \in \boldsymbol{L^2}(Q)$ is the given data for the forward problem, which can be expanded into a Fourier series, and $\boldsymbol{u}_N$ is its truncated Fourier expansion.
The residual functions depending on the Fourier coefficients are given by
\begin{align*}
 \begin{aligned}
 {\mathcal{R}_1}^c_0(\boldsymbol{\eta}_0^c,\boldsymbol{\tau}_0^c)
 &= - \text{\textbf{curl}} \, \boldsymbol{\tau}_0^c + \boldsymbol{u}^c_0, \qquad 
 {\mathcal{R}_2}^c_0(\boldsymbol{\eta}_0^c,\boldsymbol{\tau}_0^c)
 = \boldsymbol{\tau}_0^c - \nu \, \text{\textbf{curl}} \, \boldsymbol{\eta}_0^c,
 \end{aligned}
\end{align*}
and
 \begin{align*}
 {\mathcal{R}_1}_k(\boldsymbol{\eta}_k,\boldsymbol{\tau}_k)
 &= -k \omega \, \sigma \boldsymbol{\eta}_k^\perp - \text{\textbf{curl}} \, \boldsymbol{\tau}_k + \boldsymbol{u}_k, \qquad
 {\mathcal{R}_2}_k(\boldsymbol{\eta}_k,\boldsymbol{\tau}_k)
 = \boldsymbol{\tau}_k - \nu  \, \text{\textbf{curl}} \, \boldsymbol{\eta}_k.
\end{align*}

\section{A posteriori error estimates for the optimal control problem}
\label{Sec7:FunctionalAPosterioriEstimates:OptiSys}

We introduce the approximations for state $\boldsymbol{y}$ and adjoint state $\boldsymbol{p}$ denoted by $\boldsymbol{\eta}$ and $\boldsymbol{\zeta}$ from the space $\boldsymbol{H}^{\textbf{curl},1}_{0,per}(Q)$. For instance, the multiharmonic finite elements approximations \eqref{definition:MultiharmonicFEAnsatzState} and \eqref{definition:MultiharmonicFEAnsatzStateCostate} are suitable. As in the previous section, we derive an error estimate for the errors $\boldsymbol{y}-\boldsymbol{\eta}$ and $\boldsymbol{p}-\boldsymbol{\zeta}$ in $\boldsymbol{H}^{\textbf{curl},\frac{1}{2}}_{0,per}(Q)$.
The bilinear form $\mathcal{B}((\boldsymbol{y}-\boldsymbol{\eta},\boldsymbol{p}-\boldsymbol{\zeta}),(\boldsymbol{v},\boldsymbol{q}))$ defined in \eqref{definition:KKTSysSTVF} can be represented as
\begin{align}
 \label{problem:KKTSysSTVFAPost:Error}
 \begin{aligned}
  \int_{Q} \Big(&(\boldsymbol{y}-\boldsymbol{\eta})\cdot\boldsymbol{v} - \nu \, \textbf{curl} \, (\boldsymbol{p}-\boldsymbol{\zeta}) \cdot \textbf{curl} \, \boldsymbol{v}  + \sigma  \partial_t^{1/2} (\boldsymbol{p}-\boldsymbol{\zeta}) \cdot \partial_t^{1/2} \boldsymbol{v}^{\perp} \\
  &+ \nu \, \textbf{curl} \, (\boldsymbol{y}-\boldsymbol{\eta}) \cdot \textbf{curl} \, \boldsymbol{q} + \sigma \partial_t^{1/2} (\boldsymbol{y}-\boldsymbol{\eta}) \cdot \partial_t^{1/2} \boldsymbol{q}^{\perp} + \frac{1}{\alpha} \, (\boldsymbol{p}-\boldsymbol{\zeta}) \cdot\boldsymbol{q} \Big)\,d\boldsymbol{x}\,dt \\
  = \int_{Q} \Big(&\boldsymbol{y_d}\cdot\boldsymbol{v}-\boldsymbol{\eta}\cdot\boldsymbol{v} + \nu \, \textbf{curl} \, \boldsymbol{\zeta} \cdot \textbf{curl} \, \boldsymbol{v}  - \sigma  \partial_t^{1/2} \boldsymbol{\zeta} \cdot \partial_t^{1/2} \boldsymbol{v}^{\perp} \\
  &- \nu \, \textbf{curl} \, \boldsymbol{\eta} \cdot \textbf{curl} \, \boldsymbol{q} - \sigma \partial_t^{1/2} \boldsymbol{\eta} \cdot \partial_t^{1/2} \boldsymbol{q}^{\perp} - \frac{1}{\alpha} \, \boldsymbol{\zeta} \cdot\boldsymbol{q} \Big)\,d\boldsymbol{x}\,dt
 \end{aligned}
 \end{align}
for all $\boldsymbol{v}, \boldsymbol{q} \in \boldsymbol{H}^{\textbf{curl},\frac{1}{2}}_{0,per}(Q)$.
We define the right-hand side of \eqref{problem:KKTSysSTVFAPost:Error} as a linear functional of $\boldsymbol{v}$ and $\boldsymbol{q}$ as follows
\begin{align*}
 \begin{aligned}
 \mathcal{F}_{(\boldsymbol{\eta},\boldsymbol{\zeta})}(\boldsymbol{v},\boldsymbol{q})
 = \int_{Q} \Big(&\boldsymbol{y_d}\cdot\boldsymbol{v}-\boldsymbol{\eta}\cdot\boldsymbol{v} + \nu \, \textbf{curl} \, \boldsymbol{\zeta} \cdot \textbf{curl} \, \boldsymbol{v}  - \sigma  \partial_t^{1/2} \boldsymbol{\zeta} \cdot \partial_t^{1/2} \boldsymbol{v}^{\perp} \\
  &- \nu \, \textbf{curl} \, \boldsymbol{\eta} \cdot \textbf{curl} \, \boldsymbol{q} - \sigma \partial_t^{1/2} \boldsymbol{\eta} \cdot \partial_t^{1/2} \boldsymbol{q}^{\perp} - \frac{1}{\alpha} \, \boldsymbol{\zeta} \cdot\boldsymbol{q} \Big)\,d\boldsymbol{x}\,dt.
 \end{aligned}
\end{align*}
\begin{theorem}
\label{theorem:aposteriorEstimateH11/2Seminorm:OCP}
 Let $\boldsymbol{\eta}, \boldsymbol{\zeta} \in \boldsymbol{H}^{\textbf{curl},1}_{0,per}(Q)$.
 Let the bilinear form $\mathcal{B}(\cdot,\cdot)$ fulfill the inf-sup condition in \eqref{inequality:KKTSysSTVFinfsupsupsup}. The error between the exact state $(\boldsymbol{y}-\boldsymbol{\eta},\boldsymbol{p}-\boldsymbol{\zeta})$ can be estimated from above as follows
 \begin{align}
  \label{inequality:aposteriorEstimateH11/2Seminorm:OCP}
  \begin{aligned}
  \|(\boldsymbol{y}-\boldsymbol{\eta},\boldsymbol{p}-\boldsymbol{\zeta})\|_{\boldsymbol{H}^{\textbf{curl},\frac{1}{2}}}
  \leq \frac{1}{\underline{c}} &\Big( \|\mathcal{R}_1(\boldsymbol{\eta},\boldsymbol{\zeta},\boldsymbol{\tau})\|^2 +\|\mathcal{R}_3(\boldsymbol{\eta},\boldsymbol{\zeta},\boldsymbol{\rho})\|^2 \\
  &+ \|\mathcal{R}_2(\boldsymbol{\eta},\boldsymbol{\tau})\|^2 + \|\mathcal{R}_4(\boldsymbol{\zeta},\boldsymbol{\rho})\|^2\Big)^{1/2}
 \end{aligned}
 \end{align}
 with $\boldsymbol{\tau}, \boldsymbol{\rho} \in \boldsymbol{H}^{\textbf{curl},0}(Q)$ and the constant 
 $\underline{c} = (1+2\max\{\alpha,\frac{1}{\alpha}\})^{-1/2}(\min\{\frac{1}{\sqrt{\alpha}},\underline{\nu},\underline{\sigma}\} \min\{\sqrt{\alpha},\frac{1}{\sqrt{\alpha}}\})$.
 The residual functions $\mathcal{R}_1$, $\mathcal{R}_2$, $\mathcal{R}_3$ and $\mathcal{R}_4$ are defined as follows
  \begin{align}
\label{def:R1R2R3R4:1}
 \mathcal{R}_1(\boldsymbol{\eta},\boldsymbol{\zeta},\boldsymbol{\rho}) &= \sigma \partial_t \boldsymbol{\zeta} - \textbf{curl} \, \boldsymbol{\rho} + \boldsymbol{\eta} - \boldsymbol{y_d}, \quad \,\,
\mathcal{R}_2(\boldsymbol{\eta},\boldsymbol{\tau}) = \boldsymbol{\tau}
 - \nu \, \textbf{curl} \, \boldsymbol{\eta}, \\
 \label{def:R1R2R3R4:2}
\mathcal{R}_3(\boldsymbol{\eta},\boldsymbol{\zeta},\boldsymbol{\tau}) &= \sigma \partial_t \boldsymbol{\eta} + \textbf{curl} \, \boldsymbol{\tau}
 + \alpha^{-1} \boldsymbol{\zeta}, \qquad
\mathcal{R}_4(\boldsymbol{\zeta},\boldsymbol{\rho}) = \boldsymbol{\rho}
 - \nu  \, \textbf{curl} \,  \boldsymbol{\zeta}.
 \end{align}
\begin{proof}
Using identities \eqref{equation:identityEtaH11/2} and \eqref{identity:curl} for $\boldsymbol{\tau}, \boldsymbol{\rho} \in \boldsymbol{H}^{\textbf{curl},0}(Q)$, we can represent the functional $\mathcal{F}_{(\boldsymbol{\eta},\boldsymbol{\zeta})}(\boldsymbol{v},\boldsymbol{q})$ as
\begin{align*}
\mathcal{F}_{(\boldsymbol{\eta},\boldsymbol{\zeta})}(\boldsymbol{v},\boldsymbol{q}) = \int_{Q} \Big(&\boldsymbol{y_d}\cdot\boldsymbol{v}-\boldsymbol{\eta}\cdot\boldsymbol{v} + \nu \, \textbf{curl} \, \boldsymbol{\zeta} \cdot \textbf{curl} \, \boldsymbol{v}  - \sigma  \partial_t \boldsymbol{\zeta} \cdot \boldsymbol{v} + (\boldsymbol{\tau} \cdot \textbf{curl} \, \boldsymbol{q} - \textbf{curl} \, \boldsymbol{\tau} \cdot  \, \boldsymbol{q}) \\
  &- \nu \, \textbf{curl} \, \boldsymbol{\eta} \cdot \textbf{curl} \, \boldsymbol{q} - \sigma \partial_t \boldsymbol{\eta} \cdot \boldsymbol{q} - \frac{1}{\alpha} \, \boldsymbol{\zeta} \cdot\boldsymbol{q} + (-\boldsymbol{\rho} \cdot \textbf{curl} \, \boldsymbol{v} + \textbf{curl} \, \boldsymbol{\rho} \cdot \, \boldsymbol{v}) \Big)\,d\boldsymbol{x}\,dt \\
 = \int_{Q} \Big(&(\boldsymbol{y_d}-\boldsymbol{\eta} - \sigma  \partial_t \boldsymbol{\zeta} + \textbf{curl} \, \boldsymbol{\rho}) \cdot \, \boldsymbol{v} 
 + (\boldsymbol{\tau}- \nu \, \textbf{curl} \, \boldsymbol{\eta}) \cdot \textbf{curl} \, \boldsymbol{q}
  \\
  & + ( - \sigma \partial_t \boldsymbol{\eta} - \frac{1}{\alpha} \, \boldsymbol{\zeta} - \textbf{curl} \, \boldsymbol{\tau}) \cdot  \, \boldsymbol{q} 
  + (\nu \, \textbf{curl} \, \boldsymbol{\zeta} -\boldsymbol{\rho}) \cdot \textbf{curl} \, \boldsymbol{v}  \Big)\,d\boldsymbol{x}\,dt.
\end{align*}
We obtain the following estimate by applying the Cauchy-Schwarz inequality:
\begin{align*}
 \mathcal{F}_{(\boldsymbol{\eta},\boldsymbol{\zeta})}(\boldsymbol{v},\boldsymbol{q})
 \leq \,\|\mathcal{R}_1(\boldsymbol{\eta},\boldsymbol{\zeta},\boldsymbol{\rho})\| 
 \|\boldsymbol{v}\| 
 + \|\mathcal{R}_2(\boldsymbol{\eta},\boldsymbol{\tau})\| 
 \|\textbf{curl} \, \boldsymbol{q}\| 
 + \|\mathcal{R}_3(\boldsymbol{\eta},\boldsymbol{\zeta},\boldsymbol{\tau})\| 
 \|\boldsymbol{q}\| 
 + \|\mathcal{R}_4(\boldsymbol{\zeta},\boldsymbol{\rho})\| 
 \|\textbf{curl} \, \boldsymbol{v}\|, 
\end{align*}
with the residual functions $\mathcal{R}_1$, $\mathcal{R}_2$, $\mathcal{R}_3$ and $\mathcal{R}_4$ defined in \eqref{def:R1R2R3R4:1}-\eqref{def:R1R2R3R4:2}.
Together with 
\begin{align*}
 \mathcal{F}_{(\boldsymbol{\eta},\boldsymbol{\zeta})}(\boldsymbol{v},\boldsymbol{q})
 \leq \,\sqrt{\|\mathcal{R}_1(\boldsymbol{\eta},\boldsymbol{\zeta},\boldsymbol{\rho})\|^2
 + \|\mathcal{R}_2(\boldsymbol{\eta},\boldsymbol{\tau})\|^2
 + \|\mathcal{R}_3(\boldsymbol{\eta},\boldsymbol{\zeta},\boldsymbol{\tau})\|^2
 + \|\mathcal{R}_4(\boldsymbol{\zeta},\boldsymbol{\rho})\|^2}
 \|(\boldsymbol{v},\boldsymbol{q})\|_{\boldsymbol{H}^{\textbf{curl},0}}
\end{align*}
and $\|(\boldsymbol{v},\boldsymbol{q})\|_{\boldsymbol{H}^{\textbf{curl},0}} \leq \|(\boldsymbol{v},\boldsymbol{q})\|_{\boldsymbol{H}^{\textbf{curl},\frac{1}{2}}}$, applying \eqref{inequality:KKTSysSTVFinfsupsupsup} and \eqref{problem:KKTSysSTVFAPost:Error} to
\begin{align*}
\underline{c}
\|(\boldsymbol{y}-\boldsymbol{\eta},\boldsymbol{p}-\boldsymbol{\zeta})\|_{\boldsymbol{H}^{\textbf{curl},\frac{1}{2}}} &\leq 
\sup_{0 \not= (\boldsymbol{v},\boldsymbol{q}) \in (\boldsymbol{H}^{\textbf{curl},\frac{1}{2}}_{0,per}(Q))^2}    \frac{\mathcal{B}((\boldsymbol{y}-\boldsymbol{\eta},\boldsymbol{p}-\boldsymbol{\zeta}),(\boldsymbol{v},\boldsymbol{q}))}{\|(\boldsymbol{v},\boldsymbol{q})\|_{\boldsymbol{H}^{\emph{\textbf{curl}},\frac{1}{2}}}} \\
&= \sup_{0 \not= (\boldsymbol{v},\boldsymbol{q}) \in (\boldsymbol{H}^{\textbf{curl},\frac{1}{2}}_{0,per}(Q))^2}    \frac{\mathcal{F}_{(\boldsymbol{\eta},\boldsymbol{\zeta})}(\boldsymbol{v},\boldsymbol{q})}{\|(\boldsymbol{v},\boldsymbol{q})\|_{\boldsymbol{H}^{\emph{\textbf{curl}},\frac{1}{2}}}}
\end{align*}
we arrive at \eqref{inequality:aposteriorEstimateH11/2Seminorm:OCP}.
\end{proof}
 \end{theorem}
 
We derive the inf-sup and sup-sup conditions also for the seminorm providing error estimates for the seminorm of the optimality system. We can prove that
 \begin{align}
  \label{inequality:KKTSysSTVFinfsupsupsup:seminorm}
  \underline{c} |(\boldsymbol{y},\boldsymbol{p})|_{\boldsymbol{H}^{\emph{\textbf{curl}},\frac{1}{2}}} \leq
  \sup_{0 \not= (\boldsymbol{v},\boldsymbol{q}) \in (\boldsymbol{H}^{\textbf{curl},\frac{1}{2}}_{0,per}(Q))^2}
    \frac{\mathcal{B}((\boldsymbol{y},\boldsymbol{p}),(\boldsymbol{v},\boldsymbol{q}))}{
    |(\boldsymbol{v},\boldsymbol{q})|_{\boldsymbol{H}^{\emph{\textbf{curl}},\frac{1}{2}}}} \leq
  \overline{c} |(\boldsymbol{y},\boldsymbol{p})|_{\boldsymbol{H}^{\emph{\textbf{curl}},\frac{1}{2}}}
 \end{align}
 for all $(\boldsymbol{y},\boldsymbol{p}) \in (\boldsymbol{H}^{\textbf{curl},\frac{1}{2}}_{0,per}(Q))^2$, where the constants are now given as $\underline{c} = \frac{1}{\sqrt{2}}\min\{\underline{\nu},\underline{\sigma}\}\min\{\alpha,\frac{1}{\alpha}\}$ and $\overline{c} = \max\{1, 1+ {C_F^{\text{curl}}}^2\} \max\{1,\frac{1}{\alpha},\overline{\nu},\overline{\sigma}\}$. The proof follows the one of Lemma \ref{lemma:STBFinfsupsupsup:KKT}, where additionally we estimate
\begin{align*}
  \big|\mathcal{B}((\boldsymbol{y},\boldsymbol{p}),(\boldsymbol{v},\boldsymbol{q}))\big|
  \leq &\, \max\{1,\frac{1}{\alpha},\overline{\nu},\overline{\sigma}\}
  |(\boldsymbol{y},\boldsymbol{p})|_{\boldsymbol{H}^{\emph{\textbf{curl}},\frac{1}{2}}} |(\boldsymbol{v},\boldsymbol{q})|_{\boldsymbol{H}^{\emph{\textbf{curl}},\frac{1}{2}}} \\
  \leq &\, \max\{1,\frac{1}{\alpha},\overline{\nu},\overline{\sigma}\}
  \big((1+ {C_F^{\text{curl}}}^2) \|\textbf{curl} \, \boldsymbol{p}\|^2 
      + \big\|\partial^{1/2}_t \boldsymbol{p}\big\|^2 + (1+ {C_F^{\text{curl}}}^2) \|\textbf{curl} \, \boldsymbol{y}\|^2 \\
      & + \big\|\partial^{1/2}_t \boldsymbol{y}\big\|^2\big)^{1/2}
  \big((1+ {C_F^{\text{curl}}}^2) \|\textbf{curl} \, \boldsymbol{v}\|^2
      + \big\|\partial^{1/2}_t \boldsymbol{v}\big\|^2
      + (1+ {C_F^{\text{curl}}}^2) \|\textbf{curl} \, \boldsymbol{q}\|^2 \\
      & + \big\|\partial^{1/2}_t \boldsymbol{q}\big\|^2\big)^{1/2}
   \leq \, \overline{c} \, |(\boldsymbol{y},\boldsymbol{p})|_{\boldsymbol{H}^{\emph{\textbf{curl}},\frac{1}{2}}} |(\boldsymbol{v},\boldsymbol{q})|_{\boldsymbol{H}^{\emph{\textbf{curl}},\frac{1}{2}}}
\end{align*}
by applying the Friedrichs inequality \eqref{inequality:Friedrichs:FourierSpace}.
The lower bound is computed by inserting the simplified test function $ (\boldsymbol{v},\boldsymbol{q}) = (-\frac{1}{\sqrt{\alpha}} (\boldsymbol{p} + \boldsymbol{p}^\perp), \sqrt{\alpha} (\boldsymbol{y} - \boldsymbol{y}^\perp))$ yielding
\begin{align*}
|(\boldsymbol{v},\boldsymbol{q})|_{\boldsymbol{H}^{\emph{\textbf{curl}},\frac{1}{2}}} \leq \sqrt{2} \max\{\sqrt{\alpha},\frac{1}{\sqrt{\alpha}}\} |(\boldsymbol{y},\boldsymbol{p})|_{\boldsymbol{H}^{\emph{\textbf{curl}},\frac{1}{2}}}
\end{align*}
and
\begin{align*}
 \mathcal{B}((\boldsymbol{y},\boldsymbol{p}),(\boldsymbol{v},\boldsymbol{q}))
  = &\, \frac{1}{\sqrt{\alpha}} (\nu \, \textbf{curl} \, \boldsymbol{p},\textbf{curl} \, \boldsymbol{p})
  + \sqrt{\alpha} (\nu \, \textbf{curl} \, \boldsymbol{y},\textbf{curl} \, \boldsymbol{y}) \\
  &+ \frac{1}{\sqrt{\alpha}} (\sigma \partial_t^{1/2} \boldsymbol{p}, \partial_t^{1/2} \boldsymbol{p})
    + \sqrt{\alpha} (\sigma \partial_t^{1/2} \boldsymbol{y}, \partial_t^{1/2} \boldsymbol{y}) \\
  \geq &\, \frac{\underline{\nu}}{\sqrt{\alpha}} \| \textbf{curl} \, \boldsymbol{p} \|^2
  + \underline{\nu} \sqrt{\alpha} \| \textbf{curl} \, \boldsymbol{y} \|^2 
  + \frac{\underline{\sigma}}{\sqrt{\alpha}}  \big\|\partial^{1/2}_t \boldsymbol{p}\big\|^2
    + \underline{\sigma}\sqrt{\alpha} \big\|\partial^{1/2}_t \boldsymbol{y}\big\|^2.
\end{align*}
Using both estimates leads to the inf-sup condition in \eqref{inequality:KKTSysSTVFinfsupsupsup:seminorm} with constant $\underline{c}$.
Applying the Friedrichs inequality \eqref{inequality:Friedrichs:FourierSpace} to the functional $\mathcal{F}_{(\boldsymbol{\eta},\boldsymbol{\zeta})}(\boldsymbol{v},\boldsymbol{q})$ as follows
\begin{align*}
 \mathcal{F}_{(\boldsymbol{\eta},\boldsymbol{\zeta})}(\boldsymbol{v},\boldsymbol{q})
 \leq &\,{C_F^{\text{curl}}}\|\mathcal{R}_1(\boldsymbol{\eta},\boldsymbol{\zeta},\boldsymbol{\rho})\| 
 \|\textbf{curl} \,\boldsymbol{v}\| 
 + \|\mathcal{R}_2(\boldsymbol{\eta},\boldsymbol{\tau})\| 
 \|\textbf{curl} \, \boldsymbol{q}\| \\
 &+ {C_F^{\text{curl}}} \|\mathcal{R}_3(\boldsymbol{\eta},\boldsymbol{\zeta},\boldsymbol{\tau})\| 
 \|\textbf{curl} \,\boldsymbol{q}\| 
 + \|\mathcal{R}_4(\boldsymbol{\zeta},\boldsymbol{\rho})\| 
 \|\textbf{curl} \, \boldsymbol{v}\|, 
\end{align*}
leads to the majorant for the seminorm
 \begin{align}
  \label{def:PDEconstraint:majorant:seminorm:OCP}
  \mathcal{M}^\oplus_{|\cdot|}(\boldsymbol{\eta},\boldsymbol{\zeta},\boldsymbol{\tau},\boldsymbol{\rho}) = \frac{1}{\underline{c}} \big({C_F^{\text{curl}}}( \|\mathcal{R}_1(\boldsymbol{\eta},\boldsymbol{\zeta},\boldsymbol{\tau})\| +\|\mathcal{R}_3(\boldsymbol{\eta},\boldsymbol{\zeta},\boldsymbol{\rho})\|) 
  + \|\mathcal{R}_2(\boldsymbol{\eta},\boldsymbol{\tau})\| + \|\mathcal{R}_4(\boldsymbol{\zeta},\boldsymbol{\rho})\|\big).
 \end{align}
The majorant can be estimated by it's quadratic representative applying Young's inequality:
 \begin{align}
 \label{maj:quadratic:OCP}
 \begin{aligned}
  \mathcal{M}^\oplus_{|\cdot|}(\boldsymbol{\eta},\boldsymbol{\zeta},\boldsymbol{\tau},\boldsymbol{\rho})^2 \leq &\,
  \mathcal{M}^\oplus_{|\cdot|}(\beta_1, \beta_2, \beta_3; \boldsymbol{\eta},\boldsymbol{\zeta},\boldsymbol{\tau},\boldsymbol{\rho})^2 \\
  = &\,\frac{1}{\underline{c}^2} \big({C_F^{\text{curl}}}^2 (1+\beta_1)(1+\beta_2) \|\mathcal{R}_1(\boldsymbol{\eta},\boldsymbol{\zeta},\boldsymbol{\tau})\|^2 \\
  &+{C_F^{\text{curl}}}^2 \frac{(1+\beta_1)(1+\beta_3)}{\beta_1} \|\mathcal{R}_3(\boldsymbol{\eta},\boldsymbol{\zeta},\boldsymbol{\rho})\|^2 \\
  &+ \frac{(1+\beta_1)(1+\beta_2)}{\beta_2} \|\mathcal{R}_2(\boldsymbol{\eta},\boldsymbol{\tau})\|^2 +  \frac{(1+\beta_1)(1+\beta_3)}{\beta_1 \beta_3} \|\mathcal{R}_4(\boldsymbol{\zeta},\boldsymbol{\rho})\|^2\big).
  \end{aligned}
 \end{align}
with the constants $\beta_1, \beta_2, \beta_3 > 0$.
The multiharmonic representation of \eqref{def:PDEconstraint:majorant:norm:OCP} is given by
\begin{align*}
  \mathcal{M}^\oplus_{\|\cdot\|}(\boldsymbol{\eta},\boldsymbol{\zeta},\boldsymbol{\tau},\boldsymbol{\rho}) 
 = &\, \frac{1}{\underline{c}}
  \Big(T \big(\|{\mathcal{R}_1}^c_0(\boldsymbol{\eta}_0^c,\boldsymbol{\rho}_0^c)\|_{\Omega}^2 
  + \|{\mathcal{R}_2}^c_0(\boldsymbol{\eta}_0^c,\boldsymbol{\tau}_0^c)\|_{\Omega}^2 
  + \|{\mathcal{R}_3}^c_0(\boldsymbol{\zeta}_0^c,\boldsymbol{\tau}_0^c)\|_{\Omega}^2 \\
  &+ \|{\mathcal{R}_4}^c_0(\boldsymbol{\zeta}_0^c,\boldsymbol{\rho}_0^c)\|_{\Omega}^2 \big)
  + \frac{T}{2} \sum_{k=1}^N
  \big(\|{\mathcal{R}_1}_k(\boldsymbol{\eta}_k,\boldsymbol{\zeta}_k,\boldsymbol{\rho}_k)\|_{\Omega}^2 
  + \|{\mathcal{R}_2}_k(\boldsymbol{\eta}_k,\boldsymbol{\tau}_k)\|_{\Omega}^2  \\
  &+ \|{\mathcal{R}_3}_k(\boldsymbol{\eta}_k,\boldsymbol{\zeta}_k,\boldsymbol{\tau}_k)\|_{\Omega}^2 + \|{\mathcal{R}_4}_k(\boldsymbol{\zeta}_k,\boldsymbol{\rho}_k)\|_{\Omega}^2
  \big)
  + \mathcal{Q}_N \Big)^{1/2}.
 \end{align*}
The term $\mathcal{Q}_N$ is called the remainder term and is fully computable
\begin{align*}
\mathcal{Q}_N =
 \frac{T}{2} \sum_{k=N+1}^\infty \|\boldsymbol{y_d}_k\|_{\Omega}^2
 = \frac{T}{2} \sum_{k=N+1}^\infty \left(\|\boldsymbol{y_d}_k^c\|_{\Omega}^2
 + \|\boldsymbol{y_d}_k^s\|_{\Omega}^2\right)
 = \|\boldsymbol{y_d}- \boldsymbol{y_d}_N\|
\end{align*}
where $\boldsymbol{y_d} \in \boldsymbol{L^2}(Q)$ is the given desired state, which can be expanded into a Fourier series, and $\boldsymbol{y_d}_N$ is its truncated Fourier expansion.
The residual functions depending on the Fourier coefficients are given by
\begin{align*}
 \begin{aligned}
 {\mathcal{R}_1}^c_0(\boldsymbol{\eta}_0^c,\boldsymbol{\rho}_0^c)
 &= - \text{\textbf{curl}} \, \boldsymbol{\rho}_0^c + \boldsymbol{\eta}_0^c - \boldsymbol{y_d}^c_0, \qquad 
 {\mathcal{R}_2}^c_0(\boldsymbol{\eta}_0^c,\boldsymbol{\tau}_0^c)
 = \boldsymbol{\tau}_0^c - \nu \, \text{\textbf{curl}} \, \boldsymbol{\eta}_0^c, \\
 {\mathcal{R}_3}^c_0(\boldsymbol{\zeta}_0^c,\boldsymbol{\tau}_0^c)
 &= \text{\textbf{curl}} \, \boldsymbol{\tau}_0^c + \alpha^{-1} \boldsymbol{\zeta}_0^c, \qquad 
 {\mathcal{R}_4}^c_0(\boldsymbol{\zeta}_0^c,\boldsymbol{\rho}_0^c)
 = \boldsymbol{\rho}_0^c - \nu \, \text{\textbf{curl}} \, \boldsymbol{\zeta}_0^c,
 \end{aligned}
\end{align*}
and
 \begin{align*}
 {\mathcal{R}_1}_k(\boldsymbol{\eta}_k,\boldsymbol{\zeta}_k,\boldsymbol{\rho}_k)
 &= -k \omega \, \sigma \boldsymbol{\zeta}_k^\perp - \text{\textbf{curl}} \, \boldsymbol{\rho}_k + \boldsymbol{\eta}_k - \boldsymbol{y_d}_k, \qquad
 {\mathcal{R}_2}_k(\boldsymbol{\eta}_k,\boldsymbol{\tau}_k)
 = \boldsymbol{\tau}_k - \nu  \, \text{\textbf{curl}} \, \boldsymbol{\eta}_k, \\
 {\mathcal{R}_3}_k(\boldsymbol{\eta}_k,\boldsymbol{\zeta}_k,\boldsymbol{\tau}_k)
 &= -k \omega \, \sigma \boldsymbol{\eta}_k^\perp + \text{\textbf{curl}} \, \boldsymbol{\tau}_k
 + \alpha^{-1} \boldsymbol{\zeta}_k, \qquad
 {\mathcal{R}_4}_k(\boldsymbol{\zeta}_k,\boldsymbol{\rho}_k)
 = \boldsymbol{\rho}_k - \nu \, \text{\textbf{curl}} \, \boldsymbol{\zeta}_k.
\end{align*}

\section{Numerical experiments}
\label{Sec8:NumericalResults}

First numerical tests are presented here. We show the results for two numerical examples: one for the state equation and one for the corresponding optimal control problem. 
We choose the conductivity and reluctivity parameters to be $\sigma = \nu = 1$ here.
The computational domain is chosen as the unit cube $\Omega = [0,1]^3$.
We compute the discretized solution by applying the MINRES method, see \cite{PaigeSaunders:1975},
together with the preconditioner 
\begin{align}
 \label{preconditioner}
 \left( \begin{array}{cc}
  \boldsymbol{K_h}+k \omega \boldsymbol{M_{h,\sigma}} & 0 \\
  0 & \boldsymbol{K_h}+k \omega \boldsymbol{M_{h,\sigma}} \\ \end{array} \right) 
\end{align}
for the saddle point system \eqref{equation:MultiFESysBlock:forward} reformulated as
\begin{align*}
 \left( \begin{array}{cc}
  k \omega \boldsymbol{M_{h,\sigma}} & -\boldsymbol{K_h}   \\
  -\boldsymbol{K_h} & -k \omega \boldsymbol{M_{h,\sigma}} \\ \end{array} \right) 
  \left( \begin{array}{c}
     \underline{\boldsymbol{y}}_k^s \\
     \underline{\boldsymbol{y}}_k^c \end{array} \right) = \left( \begin{array}{c}
     -\underline{\boldsymbol{u}}_k^c \\
     -\underline{\boldsymbol{u}}_k^s \end{array} \right)
\end{align*}
for the forward problem.
For the optimal control problem, we use the preconditioners
\begin{align}
 \label{preconditioner:ocp}
 \left( \begin{array}{cccc}
     \boldsymbol{K_h}+k \omega \boldsymbol{M_{h,\sigma}}  &  0 & 0 & 0 \\
     0  &  \boldsymbol{K_h}+k \omega \boldsymbol{M_{h,\sigma}} & 0 & 0 \\
     0  &  0 & \alpha^{-1}(\boldsymbol{K_h}+k \omega \boldsymbol{M_{h,\sigma}}) & 0 \\
     0  &  0 & 0 & \alpha^{-1}(\boldsymbol{K_h}+k \omega \boldsymbol{M_{h,\sigma}}) \end{array} \right)
\end{align}
and
\begin{align}
 \label{preconditioner:ocp:0}
 \left( \begin{array}{cc}
  \boldsymbol{K_h}+k \omega \boldsymbol{M_{h,\sigma}} & 0 \\
  0 & \alpha^{-1}(\boldsymbol{K_h}+k \omega \boldsymbol{M_{h,\sigma}}) \\ \end{array} \right) 
\end{align}
for \eqref{equation:MultiFESysBlock} and \eqref{equation:MultiFESysBlock:forward}, the latter one for $k=0$.
Preconditioners of this type have been previously discussed in \cite{KolmbauerLanger:2012,LangerWolfmayr:2013,Wolfmayr:2014}, in the latter two for time-periodic parabolic problems.
We present the results for different Fourier modes. For evaluation of the majorant performance, we present the so called efficiency index values computed by
\begin{align*}
    I_{\text{eff}} = \frac{\mathcal{M}_{|\cdot|}^{\oplus}(\beta;\boldsymbol{\eta},\boldsymbol{\tau})^2}{\|\boldsymbol{y}-\boldsymbol{\eta}\|_{\boldsymbol{H}^{\textbf{curl},\frac{1}{2}}}}
\end{align*}
for the forward problem and 
\begin{align*}
I_{\text{eff}} = 
\frac{\mathcal{M}^\oplus_{|\cdot|}(\beta_1,\beta_2,\beta_3;\boldsymbol{\eta},\boldsymbol{\zeta},\boldsymbol{\tau},\boldsymbol{\rho})^2}{\|(\boldsymbol{y}-\boldsymbol{\eta},\boldsymbol{p}-\boldsymbol{\zeta})\|_{\boldsymbol{H}^{\textbf{curl},\frac{1}{2}}}}
\end{align*}
for the optimal control problem.
We present the results for the majorants \eqref{maj:quadratic:forward} and \eqref{maj:quadratic:OCP}, since we can use optimization with respect to the parameters $\beta$ and $\beta_1, \beta_2, \beta_3$ in order to derive better efficiency indices.

The numerical tests were computed in Matlab (MATLAB R2022a) on a computer with Intel(R) Core(TM) i5-8250U CPU @ 1.60GHz 1.80 GHz processor with 16.0 GB RAM and 512 MB system memory.
We used the \textit{Fast FEM assembly: edge elements toolbox} (\cite{Valdman:2023}) for computing the finite element discretization including mass and stiffness matrices and load vector. 

\subsection{Forward problem}
We have chosen the given data 
\begin{align*}
    \boldsymbol{u}(\boldsymbol{x},t) = (0,0,e^t (\cos{t} + (2\pi^2 +1) \sin{t}) \sin{\pi x_1} \sin{\pi x_2} )^T
\end{align*}
for which the exact solution is given by
\begin{align}
\label{numRes:exact}
    \boldsymbol{y}(\boldsymbol{x},t) = (0,0,e^t \sin{t} \sin{\pi x_1} \sin{\pi x_2} )^T.
\end{align}
In Table \ref{tab:Ex1}, we present the results for the minimization of the majorant $\mathcal{M}_{|\cdot|}^{\oplus}(\beta;\boldsymbol{\eta},\boldsymbol{\tau})^2$ with respect to $\beta$ for the mode $k=1$ including computational times in seconds $ctime$ for computing the minimization of the majorants. The corresponding error norm, which is used to compute the efficiency index of the majorant for $k=1$ is given by $1.52e+02$. The iteration stopped at iteration 8 with the stopping criterion being the value of the error between iteration steps smaller than $1e-04$.
\begin{table}[!ht]
\begin{center}
\begin{tabular}{|c|c|c|c|c|}
  \hline
   iteration & ctime & $\beta$
   & $\mathcal{M}_{|\cdot|}^{\oplus}(\beta;\boldsymbol{\eta}_1,\boldsymbol{\tau}_1)^2$ 
   & $I_{\text{eff}}$ \\
  \hline
   1 & 2.605e-03 & 1.000 & 2.246e+02 & 1.473990   \\ 
   2 & 2.119e-03 & 1.166 & 2.089e+02 & 1.370469   \\ 
   3 & 2.722e-03 & 1.275 & 2.005e+02 & 1.315729   \\   
   4 & 4.290e-03 & 1.320 & 1.975e+02 & 1.295754   \\  
   5 & 6.208e-03 & 1.332 & 1.967e+02 & 1.290793   \\    
   6 & 2.071e-03 & 1.334 & 1.966e+02 & 1.289807   \\  
   7 & 1.129e-03 & 1.335 & 1.965e+02 & 1.289624   \\  
   8 & 1.383e-03 & 1.335 & 1.965e+02 & 1.289590   \\  
  \hline
\end{tabular}
\end{center}
\caption{
The majorant $\mathcal{M}_{|\cdot|}^{\oplus}(\beta;\boldsymbol{\eta}_k,\boldsymbol{\tau}_k)^2$ for $k=1$ and the corresponding efficiency indices with respect to $\beta$ in the majorant minimization (forward problem).}
\label{tab:Ex1}
\end{table}
In Table \ref{tab:Ex1:0}, we present the results for the $k=0$ mode. The corresponding error norm, which is used to compute the efficiency index of the majorant is given by $5.18e+01$. Stopping criterion for the minimization is again being the value of the error between iteration steps smaller than $1e-04$. The results show a proper minimization with respect to $\beta$. Computational times are stay similar for all iteration steps.
\begin{table}[!ht]
\begin{center}
\begin{tabular}{|c|c|c|c|c|}
  \hline
   iteration & ctime & $\beta$
   & $\mathcal{M}_{|\cdot|}^{\oplus}(\beta;\boldsymbol{\eta}_0^c,\boldsymbol{\tau}_0^c)^2$ 
   & $I_{\text{eff}}$ \\
  \hline
   1 & 2.691e-03 & 1.000 & 1.093e+02 & 2.111 \\ 
   2 & 2.669e-03 & 1.725 & 1.015e+02 & 1.959 \\
   3 & 2.217e-03 & 1.968 & 1.011e+02 & 1.951 \\ 
  \hline
\end{tabular}
\end{center}
\caption{
The majorant $\mathcal{M}_{|\cdot|}^{\oplus}(\beta;\boldsymbol{\eta}_0^c,\boldsymbol{\tau}_0^c)^2$ for $k=0$ and the corresponding efficiency indices with respect to $\beta$ in the majorant minimization (forward problem).}
\label{tab:Ex1:0}
\end{table}

\subsection{Optimal control problem}
The given desired state is chosen as
\begin{align*}
    \boldsymbol{y_d}(\boldsymbol{x},t) = (0,0,e^t (\sin{t} + (2\pi^2 +1) ((2\pi^2+1)\sin{t}-\cos{t})) \sin{\pi x_1} \sin{\pi x_2} )^T
\end{align*}
which has the same solution \eqref{numRes:exact} for the cost parameter $\alpha=1$.
In Table \ref{tab:Ex2:0}, we present the majorant values, efficiency indices and computational times in seconds $ctime$ for different cost parameter values $\alpha$ for mode $k=0$.
In Table \ref{tab:Ex2}, we present the majorant values, efficiency indices and computational times in seconds $ctime$ for different cost parameter values $\alpha$ for mode $k=1$.
The minimization with respect to $\beta_1, \beta_2, \beta_3$ stopped usually after two iteration steps.
The results show an efficient performance for the majorant. The robustness with respect to different cost parameter values matches the results regarding the preconditioners \eqref{preconditioner:ocp:0} and \eqref{preconditioner:ocp} for applying preconditoned MINRES on the discretized systems of the optimality system. Computational times stay similar for all single computations.
\begin{table}[!ht]
\begin{center}
\begin{tabular}{|c|c|c|c|}
  \hline
   $\alpha$ & ctime  
   & ${\mathcal{M}_{|\cdot|}^{\oplus}}^2$ 
   & $I_{\text{eff}}$ \\
  \hline
   1e-04 & 3.466e-03 & 2.851e+05 & 1.788 \\
   1e-03 & 3.350e-03 & 2.854e+05 & 1.790 \\
   1e-02 & 4.157e-03 & 2.854e+05 & 1.790 \\
   1e-01 & 3.398e-03 & 2.853e+05 & 1.790 \\
   1e+00 & 4.364e-03 & 2.849e+05 & 1.787 \\
   1e+01 & 4.076e-03 & 2.805e+05 & 1.759 \\
   1e+02 & 3.156e-03 & 2.431e+05 & 1.518 \\
   1e+03 & 3.305e-03 & 2.345e+05 & 1.471 \\
   1e+04 & 4.958e-03 & 2.329e+05 & 1.461 \\
  \hline
\end{tabular}
\end{center}
\caption{
The majorant $\mathcal{M}_{|\cdot|}^{\oplus}(\beta_1,\beta_2,\beta_3;\boldsymbol{\eta}_0^c,\boldsymbol{\tau}_0^c,\boldsymbol{\zeta}_0^c,\boldsymbol{\rho}_0^c)^2$ for $k=0$ and the corresponding efficiency indices (optimal control problem).}
\label{tab:Ex2:0}
\end{table}
\begin{table}[!ht]
\begin{center}
\begin{tabular}{|c|c|c|c|}
  \hline
   $\alpha$ & ctime  
   & ${\mathcal{M}_{|\cdot|}^{\oplus}}^2$ 
   & $I_{\text{eff}}$ \\
  \hline
   1e-04 & 3.301e-03 & 4.983e+05 & 1.819 \\
   1e-03 & 4.666e-03 & 5.317e+05 & 2.373 \\
   1e-02 & 9.359e-03 & 5.317e+05 & 2.351 \\
   1e-01 & 4.302e-03 & 5.316e+05 & 1.455 \\
   1e+00 & 4.164e-03 & 5.265e+05 & 3.279 \\
   1e+01 & 3.334e-03 & 5.216e+05 & 1.798 \\
   1e+02 & 4.553e-03 & 4.838e+05 & 1.527 \\
   1e+03 & 4.768e-03 & 4.733e+05 & 1.486 \\
   1e+04 & 3.726e-03 & 4.718e+05 & 1.480 \\
  \hline
\end{tabular}
\end{center}
\caption{
The majorant $\mathcal{M}_{|\cdot|}^{\oplus}(\beta_1,\beta_2,\beta_3;\boldsymbol{\eta}_1,\boldsymbol{\tau}_1,\boldsymbol{\zeta}_1,\boldsymbol{\rho}_1)^2$ for $k=1$ and the corresponding efficiency indices (optimal control problem).}
\label{tab:Ex2}
\end{table}

\section{Conclusions and outlook}
\label{Sec9:ConclusionsOutlook}

In this work, we present the derivation of a posteriori error estimates for time-periodic eddy current problems including a standard time-periodic boundary value problem and a corresponding optimal control problem, where the forward problem is the PDE-constraint. The estimates are guaranteed, sharp, and fully computable. We discuss also a proper discretization method, the multiharmonic finite element method, for this type of problems and present first computational results.
We have presented here the derivation of the upper bounds for the forward and optimal control problems. The derivation of lower bounds would lead to a fully computable error bound from above and below. Its topic and computational experiments are part of a subsequent work of the author.

\section*{Acknowledgment}

The author gratefully acknowledges the financial support by the Regional Council of Central Finland/Council of Tampere Region and European Regional Development Fund as part of the \textit{coADDVA - ADDing VAlue by Computing in Manufacturing} project of Jamk University of Applied Sciences.

\bibliographystyle{abbrv}
\bibliography{bibliography-new}

\end{document}